\documentclass[10pt,a4paper]{amsart}
\usepackage{amssymb,amsmath}
\usepackage[american,french]{babel}
\usepackage{amsopn}
\usepackage{graphicx}
\usepackage[ec]{aeguill}
\usepackage{fancyhdr}
\usepackage[T1]{fontenc}
\title{Contraction semigroups of elliptic quadratic differential operators}
\newcommand{\rr}{\mathbb{R}}
\newcommand{\eps}{\varepsilon}
\newcommand{\nn}{\mathbb{N}}
\newcommand{\cc}{\mathbb{C}}

\newcommand{\lde}{L^2(\rr^n)}

\def\init{\setcounter{equa}{0}}
\def\inc{\stepcounter{equa}}
\def\num{\tag{\thesection.\theequa}}
\begin{document}
\newcounter{equa}
\selectlanguage{american}
\begin{center}
{\large\textbf{CONTRACTION SEMIGROUPS OF ELLIPTIC QUADRATIC DIFFERENTIAL OPERATORS}\\
\bigskip
\medskip
Karel \textsc{Pravda-Starov}\\
\bigskip
University of California, Berkeley}
\end{center}
\bigskip
\bigskip

\newtheorem{lemma}{Lemma}[section]
\newtheorem{definition}{Definition}[section]
\newtheorem{proposition}{Proposition}[section]
\newtheorem{theorem}{Theorem}[section]

\textbf{Abstract.} We study the contraction semigroups of elliptic quadratic differential operators. Elliptic quadratic differential operators are the  
non-selfadjoint operators defined in the Weyl quantization by complex-valued elliptic quadratic symbols. We establish in this paper that under the 
assumption of ellipticity, as soon as the real part of their Weyl symbols is a non-zero non-positive quadratic form, the norm of contraction semigroups generated
by these operators decays exponentially in time.

\medskip

\noindent
\textbf{Key words.} contraction semigroups, elliptic quadratic differential operators, exponential decay, pseudospectrum, subellipticity. 

\medskip

\noindent
\textbf{2000 AMS Subject Classification.} 47D06, 47F05.

\section{Miscellaneous facts about elliptic quadratic differential operators}
\init

We study in this paper the class of elliptic quadratic differential operators. It is the class of pseudodifferential operators defined in the Weyl quantization
\begin{equation}\label{3}\inc
q(x,\xi)^w u(x) =\frac{1}{(2\pi)^n}\int_{\rr^{2n}}{e^{i(x-y).\xi}q\Big(\frac{x+y}{2},\xi\Big)u(y)dyd\xi}, \num
\end{equation}
by some symbols $q(x,\xi)$, where $(x,\xi) \in \rr^{n} \times \rr^n$ and $n \in \nn^*$, which are some \textit{complex-valued elliptic quadratic forms} i.e. complex-valued quadratic forms verifying
\begin{equation}\label{3.5}\inc
(x,\xi) \in \rr^{n} \times \rr^n, \ q(x,\xi)=0 \Rightarrow (x,\xi)=(0,0). \num
\end{equation}
Since the symbols of these operators are some quadratic forms, these are only differential operators. Indeed, the quadratic symbol $x^{\alpha} \xi^{\beta}$, with 
$(\alpha,\beta) \in \nn^{2n}$ and $|\alpha+\beta| \leq 2$, is quantized in the differential operator
$$\frac{x^{\alpha}D_x^{\beta}+D_x^{\beta} x^{\alpha}}{2}, \ D_x=i^{-1}\partial_x.$$
Since their Weyl symbols are complex-valued, elliptic quadratic differential operators are a priori non-selfadjoint. 
We prove in this paper (Theorem \ref{theorem}) that under this assumption of ellipticity, as soon as the real part of their Weyl symbols is a non-zero non-positive quadratic form, the norm of contraction semigroups generated
by these operators decays exponentially in time.
Let us first begin by setting some notations and recalling known results about these operators.

Let $q$ be a complex-valued elliptic quadratic form
\begin{eqnarray*}
q : \rr_x^n \times \rr_{\xi}^n &\rightarrow& \cc\\
 (x,\xi) & \mapsto & q(x,\xi),
\end{eqnarray*}
with $n \in \nn^*$, that is a complex-valued quadratic form verifying (\ref{3.5}). The \textit{numerical range} $\Sigma(q)$ of $q$ is defined by the subset in the 
complex plane of all values taken by this symbol  
\begin{equation}\label{9}\inc
\Sigma(q)=q(\rr_x^n \times \rr_{\xi}^n), \num
\end{equation}
and the \textit{Hamilton map} $F \in M_{2n}(\cc)$ associated to the quadratic form $q$ is uniquely defined by the identity
\begin{equation}\label{10}\inc
q\big{(}(x,\xi);(y,\eta) \big{)}=\sigma \big{(}(x,\xi),F(y,\eta) \big{)}, \ (x,\xi) \in \rr^{2n},  (y,\eta) \in \rr^{2n}, \num
\end{equation}
where $q\big{(}\textrm{\textperiodcentered};\textrm{\textperiodcentered} \big{)}$ stands for the polar form associated to the quadratic form $q$ and $\sigma$ is the symplectic
form on $\rr^{2n}$,
\begin{equation}\label{11}\inc
\sigma \big{(}(x,\xi),(y,\eta) \big{)}=\xi.y-x.\eta, \ (x,\xi) \in \rr^{2n},  (y,\eta) \in \rr^{2n}. \num
\end{equation}
Let us first notice from this definition that an Hamilton map is skew-symmetric with respect to $\sigma$. This is just a consequence of the properties of skew-symmetry of the symplectic form
and symmetry of the polar form
\begin{equation}\label{12}\inc
\forall X,Y \in \rr^{2n}, \ \sigma(X,FY)=q(X;Y)=q(Y;X)=\sigma(Y,FX)=-\sigma(FX,Y).\num
\end{equation}

Under this assumption of ellipticity, the numerical range of a quadratic form can only take very particular shapes. Indeed, J. Sjöstrand proves in the lemma 3.1 in \cite{sjostrand} that
if 
\begin{equation}\label{12.1}\inc
q :  \rr_x^n \times \rr_{\xi}^n \rightarrow \cc, \ n \geq 2, \num 
\end{equation}
is a complex-valued elliptic quadratic form then there exists $z \in \cc^*$ such that  
$\textrm{Re}(z q)$ is a positive definite quadratic form. If $n=1$, the same result is fulfilled if we assume besides that $\Sigma(q) \neq \cc$. 
This result shows that the numerical range of an elliptic quadratic form can only take two shapes. The first possible shape is when $\Sigma(q)$ is equal to the whole 
complex plane. This case can only occur in dimension $n=1$. The second possible shape is when $\Sigma(q)$ is equal to a closed angular sector with a top in $0$ and an 
opening strictly lower than $\pi$.
\begin{figure}[ht]
\caption{Shape of the numerical range $\Sigma(q)$ when $\Sigma(q) \neq \cc$.}
\centerline{\includegraphics[scale=0.9]{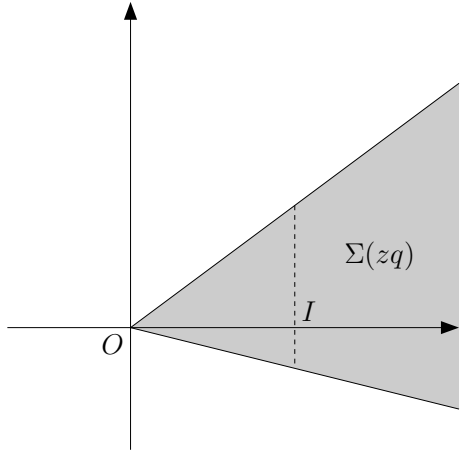}}
\end{figure}

\noindent
Indeed, if $\Sigma(q) \neq \cc$, using that the set $\Sigma(q)$ is a semi-cone  
$$t q(x,\xi)=q(\sqrt{t} x,\sqrt{t}\xi), \ t \in \rr_+, \ (x,\xi) \in \rr^{2n},$$ 
because $q$ is a quadratic form, we have  
$$\Sigma(q)=\rr_+ z^{-1}I,$$
if $z$ is the non-zero complex number mentioned above and $I$ the compact interval 
$$I=1+i \ \textrm{Im}(z q)(K),$$
where $K$ is the following compact subset of $\rr^{2n}$, 
$$\big\{(x,\xi) \in \rr^{2n} : \textrm{Re}(z q)(x,\xi)=1\big\}.$$
The compactness of $K$ is a direct consequence of the fact that $\textrm{Re}(zq)$ is a positive definite quadratic form.

Elliptic quadratic differential operators define some Fredholm operators (see Lemma~3.1 in \cite{hypoelliptic} or Theorem 3.5 in \cite{sjostrand}),
\begin{equation}\label{14}\inc
q(x,\xi)^w +z : B \rightarrow \lde, \num
\end{equation}
where $B$ is the Hilbert space 
\begin{equation}\label{14.1}\inc
\big\{ u \in \lde : x^{\alpha} D_x^{\beta} u \in \lde \ \textrm{if} \ |\alpha+\beta| \leq 2\big\}, \num
\end{equation}
with the norm  
$$\|u\|_B^2=\sum_{|\alpha+\beta| \leq 2}{\|x^{\alpha} D_x^{\beta} u\|_{\lde}^2}.$$
The Fredholm index of the operator $q(x,\xi)^w+z$ is independent of $z$ and is equal to~$0$ if $n \geq 2$. In the case where $n=1$, this index can take the values $-2$, $0$ or $2$. More 
precisely, this index is always equal to $0$ if $\Sigma(q) \neq \cc$.

If $\Sigma(q) \neq \cc$, J.~Sjöstrand proves in the theorem 3.5 in \cite{sjostrand} (see also Lemma $3.2$ and Theorem~3.3 
in~\cite{hypoelliptic}) that the spectrum of an elliptic quadratic differential operator  
$$q(x,\xi)^w : B \rightarrow \lde,$$ 
is only composed of eigenvalues with finite multiplicity 
\begin{equation}\label{15}\inc
\sigma\big{(}q(x,\xi)^w\big{)}=\Big\{ \sum_{\substack{\lambda \in \sigma(F), \\  -i \lambda \in \Sigma(q) \setminus \{0\}}}
{\big{(}r_{\lambda}+2 k_{\lambda}
\big{)}(-i\lambda) : k_{\lambda} \in \nn} 
\Big\}, \num
\end{equation}
where $F$ is the Hamilton map associated to the quadratic form $q$ and $r_{\lambda}$ is the dimension of the space of generalized eigenvectors of $F$ in $\cc^{2n}$
belonging to the eigenvalue $\lambda \in \cc$. Let us notice that the spectra of these operators are always included in the numerical range of their Weyl symbols.

\section{Contraction semigroups of elliptic quadratic differential operators}
\init

Let us consider a complex-valued elliptic quadratic form
\begin{eqnarray*}
q : \rr_x^n \times \rr_{\xi}^n &\rightarrow& \cc\\
 (x,\xi) & \mapsto & q(x,\xi),
\end{eqnarray*}
with $n \in \nn^*$, verifying 
\begin{equation}\label{inf1}\inc
\textrm{Re }q(x,\xi) \leq 0, \ (x,\xi) \in \rr^{2n}. \num
\end{equation}
Since the quadratic form $\textrm{Re } q$ is non-positive, we deduce from the theorem 21.5.3 in~~\cite{hormander} that there exists a real linear symplectic transformation $\chi$ of 
$\rr^{2n}$ such that 
\begin{equation}\label{inf2}\inc
-\textrm{Re}(q \circ \chi)(x,\xi)=\sum_{j=1}^{k}{\lambda_j(\xi_j^2+x_j^2)}
+\sum_{j=k+1}^{k+l}{x_j^2}, \num
\end{equation}
with $k,l \in \nn$ and $\lambda_j>0$ for all $j=1,...,k$. Then, it follows from 
the symplectic invariance of the Weyl quantization (Theorem 18.5.9 in 
\cite{hormander}) that we can find a metaplectic operator $U$,
which is a unitary transformation on $\lde$ and an automorphism of the spaces $\mathcal{S}(\rr^n)$ and $B$, where $B$ is the Hilbert space defined in (\ref{14.1}) such that  
\begin{equation}\label{inf3}\inc
\textrm{Re } q(x,\xi)^w=-U^{-1} \Big(\sum_{j=1}^{k}{\lambda_j(D_{x_j}^2+x_j^2)}
+\sum_{j=k+1}^{k+l}{x_j^2} \Big) U. \num
\end{equation}
By using that $U$ is a unitary operator on $L^2(\rr^{n})$, we obtain that for all $\lambda \in \rr$ and $u \in \mathcal{S}(\rr^n)$,
\inc
\begin{multline*}\label{inf3.5}
\textrm{Re}(\lambda u-q(x,\xi)^w u,u)_{L^2}=  \lambda \|u\|_{L^2}^2-(\textrm{Re }q(x,\xi)^wu,u)_{L^2} \\
=  \lambda \|u\|_{L^2}^2 +\sum_{j=1}^{k}{\lambda_j (\|D_{x_j} U u\|_{L^2}^2+\|x_j U u\|_{L^2}^2 )}+\sum_{j=k+1}^{k+l}{\|x_j U u\|_{L^2}^2} 
\geq  \lambda \|u\|_{L^2}^2, \num
\end{multline*} 
because the operator $i \textrm{Im }q(x,\xi)^w$ is formally skew-selfadjoint. 
By using the Cauchy-Schwarz inequality and the density of the Schwartz space $\mathcal{S}(\rr^n)$ in $B$, we deduce from (\ref{inf3.5}) that
$$\forall \lambda >0, \forall m \in \nn^*, \ \big\|\big(q(x,\xi)^w-\lambda\big)^{-m}\big\| \leq \lambda^{-m},$$ 
because we get from (\ref{15}) and (\ref{inf1}) that 
$$\sigma\big(q(x,\xi)^w\big) \subset \{z \in \cc : \textrm{Re } z \leq 0\}.$$
This induces that the elliptic quadratic differential operator $q(x,\xi)^w$ generates a contraction semigroup 
$$e^{t q(x,\xi)^w}, \  t \geq 0.$$
Indeed, this fact is a direct consequence of the following classical result, which can be found for example in the theorem~2.21 in \cite{davies}: let $Q$ be a closed 
unbounded linear operator with a dense domain on 
$L^2(\rr^n)$, $Q$ is the generator of a one-parameter semigroup $e^{tQ}$ verifying
\begin{equation}\label{inf4}\inc
\exists M>0, \exists a \in \rr, \forall t \geq 0, \ \|e^{tQ}\|_{\mathcal{L}(L^2)} \leq M e^{-a t}, \num
\end{equation}
if and only if 
\begin{equation}\label{inf5}\inc
\forall \lambda>-a, \forall m \in \nn^*, \ \lambda \not\in \sigma(Q) \textrm{ and } \|(Q-\lambda)^{-m}\| \leq M(\lambda+a)^{-m}. \num
\end{equation}

\medskip

\noindent
\textit{Remark.} As mentioned in \cite{mehler} (p. 426), the assumption of ellipticity for a quadratic symbol $q$ verifying (\ref{inf1}) is not necessary to obtain 
that the operator defined in the Weyl quantization generates a contraction semigroup.

\subsection{Statement of the main result}

The main result contained in this paper is the following.

\medskip

\begin{theorem}\label{theorem}
If $q : \rr_x^n \times \rr_{\xi}^n \rightarrow \cc$ is a complex-valued elliptic quadratic form with $n \in \nn^*$ such that  
\begin{equation}\label{inf6}\inc
\emph{\textrm{Re }}q \leq 0, \textrm{ and } \exists (x_0,\xi_0) \in \rr^{2n}, \ \emph{\textrm{Re }}q(x_0,\xi_0) \neq 0, \num
\end{equation}
then the elliptic quadratic differential operator $q(x,\xi)^w$ generates a contraction semigroup with a norm decaying exponentially in time
\begin{equation}\label{inf6.5}\inc
\exists M>0, \exists a >0, \forall t \geq 0, \ \|e^{t q(x,\xi)^w}\|_{\mathcal{L}(L^2)} \leq M e^{-a t}. \num
\end{equation}
\end{theorem}

\medskip

\subsection{Some remarks about this result}

Let us first notice that the assumption of ellipticity is essential. Indeed, let us consider the operator defined in the Weyl quantization by the following 
quadratic form
$$q(x,\xi)=-x^2.$$
This operator is just the operator of multiplication by $-x^2$, which generates the contraction semigroup
$$e^{t q(x,\xi)^w}u=e^{-t x^2}u, \ t \geq 0, \ u \in L^2(\rr^n).$$  
We can easily check in this case that for all $t \geq 0$,
$$\|e^{t q(x,\xi)^w}\|_{\mathcal{L}(L^2)}=1.$$
The second assumption in (\ref{inf6}), 
$$\exists (x_0,\xi_0) \in \rr^{2n}, \ \textrm{Re }q(x_0,\xi_0) \neq 0,$$
is also important. Indeed, if the quadratic form $\textrm{Re }q$ is identically equal to zero, the two elliptic quadratic differential operators $q(x,\xi)^w$ and $-q(x,\xi)^w$ generate
some contraction semigroups verifying
$$\big(e^{t q(x,\xi)^w}\big)^{-1}=e^{t (-q(x,\xi)^w)}, \ t \geq 0.$$
This induces that we also have for all $t \geq 0$,  
$$\|e^{t q(x,\xi)^w}\|_{\mathcal{L}(L^2)}=1,$$ 
in this case.

\subsection{Structure of the proof}
We shall study two cases in the proof of the theorem~\ref{theorem}. The first case is the case where the Hamilton map $\textrm{Re }F$ associated to the real 
part of the operator's Weyl symbol $\textrm{Re }q$ is \textit{non-nilpotent}. This case corresponds to the case where the integer $k$ appearing in (\ref{inf2}) is non-zero. In this case, 
the assumption 
of ellipticity will not be needed and we shall get straightforwardly the result of the theorem \ref{theorem} from (\ref{inf3.5}).

The most interesting case is the case where the Hamilton map $\textrm{Re }F$
is \textit{nilpotent}. The assumption of ellipticity is then essential. In this case, we shall see that in a suitable new system of symplectic coordinates, we can tensor 
the variables 
$$(x;\xi)=(x',x'';\xi',\xi'') \in \rr^{2n},$$ 
to write the Weyl symbol of the studied operator
$$q(x,\xi)=\tilde{q}(x',\xi')+i \eps \sum_{j=1}^{n''}{\mu_j(\xi_j''^2+x_j''^2)},$$
with $\mu_j>0$ for all $j=1,...,n''$, $\eps \in \{\pm 1\}$ and $n=n'+n''$, as a sum of a complex-valued elliptic quadratic form $\tilde{q}$ on $\rr^{2n'}$ verifying
$$\textrm{Re }\tilde{q} \leq 0,$$
such that its associated Hamilton map $\tilde{F}$ has no real eigenvalue, 
and the symbol of an harmonic oscillator in dimension $n''$. Thanks to this tensorization, it will be sufficient to prove the result of exponential decay in time for the norm of the contraction semigroup
generated by the elliptic quadratic differential operator $\tilde{q}(x',\xi')^w$.

In a second step, we shall see that the properties about its Hamilton map will induce that the symbol 
$\tilde{q}(x',\xi')$ is of \textit{finite} order in every non-zero point of its numerical range.  By using some results proved in \cite{karel} about the pseudospectrum 
of elliptic quadratic differential operators\footnote{These results are linked to some properties of subellipticity.}, we will deduce from this fact some estimates for the resolvent's norm of 
the operator $\tilde{q}(x',\xi')^w$ inside some particular regions of its numerical range $\Sigma(\tilde{q})$. These estimates will allow to give an integral formula for the contraction semigroup 
generated by the operator $\tilde{q}(x',\xi')^w$, which will permit to prove the result of exponential decay in time for its norm.


\section{Proof of the theorem \ref{theorem}}
\init

Let us denote by $F$ the Hamilton map associated to the quadratic form $q$. It follows from (\ref{10}) that its real part $\textrm{Re } F$ and its imaginary part
$\textrm{Im }F$ are respectively the Hamilton maps associated to the quadratic forms $\textrm{Re } q$ and $\textrm{Im }q$. We deduce from (\ref{10}) and (\ref{inf2}) that 
we have in the system of symplectic coordinates $(y_1,\eta_1,...,y_n,\eta_n)$, where $(y,\eta)=\chi(x,\xi)$,   
\begin{equation}\label{2.3.88}\inc
\textrm{Re } F=M \left(
  \begin{array}{ccccccc}
  J_1 & 0 & \ldots & 0 & \ldots &0\\
  0& \ddots & \ddots &\vdots &0 & \vdots \\
\vdots & \ddots & \ddots & 0 &\ldots & 0\\
0 & \ldots & 0 & J_{k+l} & 0 & 0\\
0 & \ldots & \ldots & 0 &\ldots &0 \\
0 & \ldots & \ldots & \ldots &\ldots &0 \\
0 & \ldots & \ldots & 0 &\ldots &0 \\
  \end{array}
\right) M^{-1}, \num
\end{equation}
with
\begin{equation}\label{inf7}\inc
J_j=\left(
  \begin{array}{cc}
  0 & -\lambda_j \\
  \lambda_j & 0 \\
  \end{array}
\right) \textrm{ if } 1 \leq j \leq k, \textrm{ and } 
J_j=\left(
  \begin{array}{cc}
  0 & 0 \\
  1& 0 \\
  \end{array}
\right) \textrm{ if } k+1 \leq j \leq k+l, \num
\end{equation}
if $M$ stands for the matrix of the real linear symplectic transformation $\chi$.
We deduce from (\ref{inf2}), (\ref{2.3.88}) and (\ref{inf7}) that the Hamilton map $\textrm{Re } F$ is \textit{nilpotent} if and only if $k=0$. In this case $k=0$, the index of nilpotence
is equal to 2. Let us begin by first considering the non-nilpotent case.

\subsection{The non-nilpotent case} Let us consider the case where $k \geq 1$. We deduce from (\ref{inf3.5}) that for all $\lambda \in \rr$ and $u \in \mathcal{S}(\rr^n)$,
\inc
\begin{align*}\label{inf8}
& \ \textrm{Re}(\lambda u-q(x,\xi)^w u,u)_{L^2}=  \lambda \|u\|_{L^2}^2-(\textrm{Re }q(x,\xi)^wu,u)_{L^2} \num \\
=  & \ \lambda \|u\|_{L^2}^2 +\sum_{j=1}^{k}{\lambda_j (\|D_{x_j} U u\|_{L^2}^2+\|x_j U u\|_{L^2}^2 )}+\sum_{j=k+1}^{k+l}{\|x_j U u\|_{L^2}^2} \\
\geq & \  \lambda \|u\|_{L^2}^2+\lambda_1 (\|D_{x_1} U u\|_{L^2}^2+\|x_1 U u\|_{L^2}^2), 
\end{align*}
because $\lambda_j>0$ for all $j=1,...,k$. 
Since we get from the Cauchy-Schwarz inequality
\begin{multline*}
\|D_{x_1} U u\|_{L^2}^2+\|x_1 U u\|_{L^2}^2 \geq 2\|D_{x_1} U u\|_{L^2}\|x_1 U u\|_{L^2} \\ \geq 2\textrm{Re}(D_{x_1}Uu,ix_1 U u)_{L^2}
= ([D_{x_1},ix_1]Uu,U u)_{L^2}=\|Uu\|_{L^2}^2=\|u\|_{L^2}^2, 
\end{multline*}
where $[D_{x_1},ix_1]$ stands for the commutator of the operators $D_{x_1}$ and $ix_1$, because $U$ is unitary on $L^2(\rr^n)$, we deduce from (\ref{inf8}) that 
for all $\lambda \in \rr$ and $u \in \mathcal{S}(\rr^n)$,
\begin{equation}\label{inf9}\inc
\|\lambda u-q(x,\xi)^w u\|_{L^2} \geq (\lambda+\lambda_1) \|u\|_{L^2}. \num 
\end{equation} 
Since from (\ref{15}), the spectrum of the operator $q(x,\xi)^w$ is only composed of eigenvalues and that the Schwartz space $\mathcal{S}(\rr^n)$ is dense in $B$, we get from
(\ref{inf9}) that 
$$\forall \lambda > - \lambda_1, \forall m \in \nn^*, \ \|(q(x,\xi)^w-\lambda)^{-m}\| \leq (\lambda+\lambda_1)^{-m},$$
which induces in view of (\ref{inf4}) and (\ref{inf5}) that for all $t \geq 0$,
$$\|e^{t q(x,\xi)^w}\|_{\mathcal{L}(L^2)} \leq  e^{-\lambda_1 t}.$$
This proves the theorem \ref{theorem} in the non-nilpotent case.

\subsection{The nilpotent case}
Let us now consider the nilpotent case, which is the most interesting one. The Hamilton map $\textrm{Re }F$ is then nilpotent of index 2, 
\begin{equation}\label{inf9.5}\inc
(\textrm{Re }F)^2=0.\num
\end{equation}
We deduce from (\ref{inf2}), (\ref{inf6}), (\ref{2.3.88}) and (\ref{inf7}) 
that in this case 
\begin{equation}\label{inf10}\inc
k=0 \textrm{ and } l \geq 1,\num
\end{equation}
because the quadratic form $\textrm{Re }q$ is non identically equal to zero.

\subsubsection{Tensorization of variables}\label{s3.2.1}
Let us first notice that we can freely choose the symplectic coordinates, in which we express the symbol $q$ of the studied operator. Indeed, by using 
the symplectic invariance of the Weyl quantization given by the theorem~18.5.9 in 
\cite{hormander},  we can find a metaplectic operator $U$, which is a unitary transformation on $\lde$ and an automorphism of the spaces $\mathcal{S}(\rr^n)$ and $B$ such that  
\begin{equation}\label{infa3}\inc
(q \circ \chi)(x,\xi)^w=U^{-1} q(x,\xi)^w U, \num
\end{equation}
if $\chi$ is a real linear symplectic transformation of $\rr^{2n}$.
Since (\ref{infa3}) induces that 
$$e^{t(q \circ \chi)(x,\xi)^w}=U^{-1} e^{t q(x,\xi)^w} U, \ t \geq 0,$$
at the level of the generated semigroups, we actually get that for all $t \geq 0$,
$$\|e^{t(q \circ \chi)(x,\xi)^w}\|_{\mathcal{L}(L^2)}=\|e^{t q(x,\xi)^w}\|_{\mathcal{L}(L^2)},$$
because $U$ is a unitary operator on $L^2(\rr^n)$. This last identity justifies this first remark.
 
Let us now assume that we can find a decomposition of the phase space 
$$T^* \rr^n=\rr_x^n \times \rr_{\xi}^n,$$ in a direct sum of symplectic sub-vector spaces 
$$S_j \subset T^* \rr^n, \ 1 \leq j \leq N,$$  
which are symplectically orthogonal and stable by the maps $\textrm{Re } F$ and $\textrm{Im } F$,
\begin{equation}\label{2.3.89}\inc
T^*\rr^n= \bigoplus_{j=1}^N{}^{\sigma \perp  }S_j, \  (\textrm{Re } F)S_j \subset S_j \ \textrm{and } 
(\textrm{Im } F)S_j \subset S_j, \num
\end{equation}
for all $j \in \{1,...,N\}$. We consider a symplectic basis $(e_{1,j},\eps_{1,j},...,e_{N_j,j},\eps_{N_j,j})$ of the symplectic space $S_j$. By collecting these previous basis, we get a 
symplectic basis of $T^* \rr^n$. By using the stability and the orthogonality properties of the symplectic spaces $S_j$, we obtain the following decomposition in this new system of symplectic
coordinates
\begin{align*}\label{2.3.90}\inc
& \ q(x,\xi)\\ =
& \  
\sigma\Big{(} \sum_{\substack{1 \leq j \leq N, \\ 1 \leq k \leq N_j}}{(x_{k,j} e_{k,j}+\xi_{k,j} \eps_{k,j})}, F\big{(} 
\sum_{\substack{1 \leq j \leq N,\\  1 \leq k \leq N_j}}{(x_{k,j} e_{k,j}+\xi_{k,j} \eps_{k,j})}\big)\Big)\\
= & \ 
\sum_{1 \leq j \leq N}\sigma \Big{(}\sum_{1 \leq k \leq N_j}{(x_{k,j} e_{k,j}+\xi_{k,j} \eps_{k,j})},F\big{(} 
\sum_{1 \leq k \leq N_j}{(x_{k,j} e_{k,j}+\xi_{k,j} \eps_{k,j})}\big)\Big)
\\ = & \  \sum_{1 \leq j \leq N} q_j(x_{1,j},...,x_{N_j,j},\xi_{1,j}, ...,\xi_{N_j,j}), \num
\end{align*}
with 
\inc
\begin{multline*}\label{2.3.91}
q_j(x_{1,j},...,x_{N_j,j},\xi_{1,j}, ...,\xi_{N_j,j})=\\q(0,...,0,x_{1,j},...,x_{N_j,j},0,...,0,\xi_{1,j}, ...,\xi_{N_j,j},0,...,0), \num
\end{multline*}
for all $1 \leq j \leq N$. Since from (\ref{inf6}), (\ref{2.3.90}) and (\ref{2.3.91}), $q_j$ is a complex-valued elliptic quadratic form on $\rr^{2N_j}$ verifying 
$$\textrm{Re }q_j \leq 0,$$  
the operator $q_j(x_{1,j},...,x_{N_j,j},\xi_{1,j}, ...,\xi_{N_j,j})^w$ generates of a contraction semigroup 
$$e^{t q_j(x_{1,j},...,x_{N_j,j},\xi_{1,j},...,\xi_{N_j,j})^w}, \ t \geq 0,$$ 
such that for all $t \geq 0$,
\begin{equation}\label{inf11}\inc
e^{t q(x,\xi)^w}=\prod_{j=1}^N{e^{t q_j(x_{1,j},...,x_{N_j,j},\xi_{1,j},...,\xi_{N_j,j})^w}}. \num
\end{equation}
We deduce from (\ref{inf11}) that it is sufficient for proving (\ref{inf6.5}) to find an integer $1 \leq j_0 \leq N$ such that 
$$\exists M>0, \exists a>0, \forall t \geq 0, \ \|e^{t q_{j_0}^w}\|_{\mathcal{L}(L^2(\rr^{N_{j_0}}))} \leq M e^{- a t}.$$ 
This second remark will allow us to reduce in the following our study to the case where the Hamilton map $F$ associated to the quadratic form $q$ has no real eigenvalue. 
Indeed, it will be a consequence of the following lemma.

\medskip

\begin{lemma}\label{l1}
If $\lambda$ is a real eigenvalue of the Hamilton map $F$ associated to a complex-valued elliptic quadratic form $q$ verifying \emph{(\ref{inf6})} then this eigenvalue $\lambda$ is 
necessarily non-zero 
$$\lambda \neq 0,$$
and the complex sub-vector space $\emph{\textrm{Ker}}(F+\lambda)$ is the sub-vector space complex conjugate of $\emph{\textrm{Ker}}(F-\lambda)$.
If we consider 
$$S_{\lambda}=\big{(}\emph{\textrm{Ker}}(F-\lambda) \oplus  \emph{\textrm{Ker}}(F+\lambda)\big{)} \cap T^* \rr^n,$$  
then the vector space $\emph{\textrm{Ker}}(F-\lambda) \oplus  \emph{\textrm{Ker}}(F+\lambda)$ is a complexification of the space
$S_{\lambda}$. Moreover, $S_{\lambda}$ has a structure of symplectic vector space, which is stable by the Hamilton map $\emph{\textrm{Im }} F$ and verifies 
$$\emph{\textrm{Re }} F \ S_{\lambda}=\{0\}.$$
Moreover, if $\mu$ is another real eigenvalue of $F$ distinct from $\lambda$ and $-\lambda$, then we have the following symplectically orthogonal direct sum
$$S_{\lambda} \oplus^{\sigma \perp  } S_{\mu} \subset T^* \rr^n,$$
where
$$S_{\mu}=\big{(}\emph{\textrm{Ker}}(F-\mu) \oplus  \emph{\textrm{Ker}}(F+\mu)\big{)} \cap T^* \rr^n.$$
\end{lemma}

\medskip
\noindent
\textit{Proof of Lemma \ref{l1}.} Let $\lambda$ be a real eigenvalue of $F$. The proposition $4.4$ in \cite{mehler} 
shows that the sub-vector space of $\cc^{2n}$, 
$$\textrm{Ker}(F+\lambda),$$
is the complex conjugate of the complex sub-vector space $\textrm{Ker}(F-\lambda).$
Moreover, we also get from this proposition that 
$$\textrm{Re }  F \ \textrm{Ker}(F \pm \lambda)=\{0\}$$
and that the vector space $\textrm{Ker}(F-\lambda) \oplus  \textrm{Ker}(F+\lambda)$, if $\lambda \neq 0$, and $\textrm{Ker}\ F$ are the complexifications of their intersections with $T^* \rr^n$.
If $\lambda=0$, then we could find $X_0=X_1+i X_2 \neq 0$ such that $X_1, X_2 \in T^* \rr^n$ and 
$X_0, X_1, X_2 \in \textrm{Ker } F$. Let us consider $j \in \{1,2\}$ such that $X_j \neq 0$. Thus, we would get that 
$$q(X_j)=\sigma(X_j, F X_j)=0,$$
which would contradict the ellipticity of the quadratic form $q$. This induces that the eigenvalue $\lambda$ is non-zero. The fact that  
\begin{equation}\label{2.3.93}\inc
\textrm{Re } F \ \textrm{Ker}(F \pm \lambda)=\{0\}, \num
\end{equation}
implies that 
$$\textrm{Re } F \ S_{\lambda}=\{0\}.$$
Moreover, since the space $\textrm{Ker}(F-\lambda) \oplus  \textrm{Ker}(F+\lambda)$ is stable by $F$, it follows from (\ref{2.3.93}) that it is also stable by $\textrm{Im } F$.
The stability of the sub-space $S_{\lambda}$ by the matrix $\textrm{Im } F$ then directly follows from the fact that 
$\textrm{Im } F \in M_{2n}(\rr)$. Let us now consider $X_0 \in S_{\lambda}$ such that for all $Y \in S_{\lambda}$, 
$$\sigma(X_0,Y)=0.$$ 
It follows that for all $Y$ and  $Z$ in $S_{\lambda}$,  
$$\sigma(X_0,Y+iZ)=0,$$
which induces that 
$$\forall X \in \textrm{Ker}(F-\lambda) \oplus  \textrm{Ker}(F+\lambda), \ \sigma(X_0,X)=0,$$
because $ \textrm{Ker}(F-\lambda) \oplus  \textrm{Ker}(F+\lambda)$ is a complexification of $S_{\lambda}$.
On the other hand, since $X_0 \in S_{\lambda}$, we have $F X_0 \in \textrm{Ker}(F-\lambda) \oplus  \textrm{Ker}(F+\lambda)$, which implies that  
$$q(X_0)=\sigma(X_0,F X_0)=0.$$
We can then deduce from the ellipticity of the quadratic form $q$ that $X_0=0$, which proves the symplectic structure of the space $S_{\lambda}$.

Let us now assume  
that there exists a real eigenvalue $\mu$ of $F$ distinct from $\lambda$ and $-\lambda$. We already know that this eigenvalue $\mu$ is necessarily non-zero. Let  
$$X \in \textrm{Ker}(F-\eps_1 \lambda) \textrm{ and } Y \in \textrm{Ker}(F-\eps_2 \mu),$$ 
with $\eps_1, \eps_2 \in \{\pm 1\}$, we have 
$$\sigma(X,Y)=\sigma(X,\eps_2^{-1} \mu^{-1} FY)=-\frac{1}{\eps_2 \mu}\sigma(FX,Y)=-\frac{\eps_1 \lambda}{\eps_2 \mu} \sigma(X,Y),$$
according to the skew-symmetry of the Hamilton map $F$ with respect to $\sigma$. Since 
$$\left|\frac{\eps_1 \lambda}{\eps_2 \mu} \right| \neq 1,$$
because $\lambda$ and $\mu$ are some real numbers such that $\mu \not\in \{\lambda,-\lambda\}$, we finally deduce that  
$$\sigma(X,Y) =0,$$
which proves that the sub-vector spaces $S_{\lambda}$ and $S_{\mu}$ are symplectically orthogonal. $\Box$

\bigskip
Let us consider the real eigenvalues of the Hamilton map $F$. In view of the previous lemma, we can write the two by two different elements of this set  
$$\lambda_1, ..., \lambda_r, -\lambda_1, ...,-\lambda_r.$$
By using the notations introduced above, let us consider the following space
\begin{equation}\label{2.3.94}\inc
S=S_{\lambda_1} \oplus^{\sigma \perp  } ... \oplus^{\sigma \perp  } S_{\lambda_r} \subset T^* \rr^n, \num
\end{equation}
which is actually a direct sum because of the symplectic structure of sub-vector spaces $S_{\lambda_j}$ and their properties of two by two symplectic 
orthogonality. The lemma~\ref{l1} shows that the symplectic vector space $S$ is stable by the Hamilton maps $\textrm{Re } F$ and
$\textrm{Im } F$. Let us denote by $S^{\sigma \perp}$ its symplectic orthogonal vector space in $T^* \rr^n$. This space $S^{\sigma \perp}$
is then a symplectic vector space, which is also stable by $\textrm{Re } F$ and $\textrm{Im } F$. Indeed, let us consider 
$X \in S^{\sigma \perp} \subset T^* \rr^n$. 
Since from the skew-symmetry of Hamilton maps with respect to $\sigma$, we get that for all $Y \in S \subset T^* \rr^n$,  
$$\sigma(Y,\textrm{Re } F X)+i \sigma(Y,\textrm{Im } F X)=-\sigma(\textrm{Re } F Y,X)-i \sigma(\textrm{Im } F Y,X)=0,$$
because $\textrm{Re } F Y \in S$ and $\textrm{Im } F Y \in S$. 
It follows that for all $Y \in S$, 
$$\sigma(Y,\textrm{Re } F X)=\sigma(Y,\textrm{Im } F X)=0,$$
which proves that $\textrm{Re } F X \in S^{\sigma \perp  }$ and $\textrm{Im } F X \in S^{\sigma \perp  }$. In view of our construction (\ref{2.3.94}), 
the map $\textrm{Re }F|_{S^{\sigma \perp}}+i\textrm{Im }F|_{S^{\sigma \perp}}$ has no real eigenvalue. 
If we now consider some symplectic coordinates $(x,\xi)$,   
$$x=(x';x'') \in \rr^{n'} \times \rr^{n''}, \ \xi=(\xi';\xi'') \in \rr^{n'} \times \rr^{n''},$$ 
such that $(x',\xi')$ and $(x'',\xi'')$ are respectively some symplectic coordinates in $S$ and $S^{\sigma \perp  }$, we deduce from the symplectic orthogonality 
of the spaces $S$ and $S^{\sigma \perp  }$, their stabilities by Hamilton maps $\textrm{Re }F$, $\textrm{Im }F$ and the previous results that we can 
tensorize the variables in the symbol $q$,
\begin{equation}\label{inf11.5}\inc
q(x,\xi)=\tilde{q}_1(x',\xi')+\tilde{q}_2(x'',\xi''), \num
\end{equation}
with
\begin{equation}\label{inf12}\inc
\tilde{q}_1(x',\xi')=\sigma_{n'}\big((x',\xi'),(\textrm{Re }F|_S+i \textrm{Im }F|_S)(x',\xi')\big) \num
\end{equation}
and
\begin{equation}\label{inf13}\inc
\tilde{q}_2(x'',\xi'')=\sigma_{n''} \big((x'',\xi''),(\textrm{Re }F|_{S^{\sigma \perp }}+i\textrm{Im }F|_{S^{\sigma \perp }})(x'',\xi'')\big), \num
\end{equation}
if $\sigma_{n'}$ and $\sigma_{n''}$ stand respectively for the canonical symplectic forms on $\rr^{2n'}$ and $\rr^{2n''}$. We deduce from (\ref{inf6}),
(\ref{inf12}) and (\ref{inf13}) that $\tilde{q}_1$ and $\tilde{q}_2$ are some complex-valued elliptic quadratic forms respectively on $\rr^{2n'}$ and $\rr^{2n''}$ verifying
$$\textrm{Re }\tilde{q}_j \leq 0, \ j=1,2.$$
Since from (\ref{2.3.94}) and the lemma \ref{l1}, 
$$\textrm{Re } F|_{S}=0,$$
we deduce from (\ref{inf6}), (\ref{inf11.5}), (\ref{inf12}) and (\ref{inf13}) that there exists $(x_0'',\xi_0'') \in \rr^{2n''}$ such that 
$$\textrm{Re }\tilde{q}_2(x_0'',\xi_0'') \neq 0.$$
Let us also notice that (\ref{inf9.5}) implies that
$$(\textrm{Re }F|_{S^{\sigma \perp }})^2=0.$$ 
It follows in view of the reasoning mentioned above that we can reduce our study in the nilpotent case to the case where the Hamilton map $F$ associated to the complex-valued 
quadratic form $q$ defined in (\ref{inf6}) has \textit{no real eigenvalue}. In all the following, we make this additional licit assumption.

\subsubsection{Some algebraic properties}\label{s3.2.2}
Let us begin this paragraph by proving the following lemma.

\medskip

\begin{lemma}\label{ll2}
If $q_1$ and $q_2$ are two complex-valued quadratic forms on $\rr^{2n}$, then the Hamilton map associated to the complex-valued quadratic form defined by the Poisson 
bracket 
$$\{q_1,q_2\}=\frac{\partial q_1}{\partial \xi}.\frac{\partial q_2}{\partial x}-\frac{\partial q_1}{\partial x}.\frac{\partial q_2}{\partial \xi},$$
is $-2[F_1,F_2]$ if $[F_1,F_2]$ stands for the commutator of $F_1$ and $F_2$ the Hamilton maps of $q_1$ and $q_2$.  
\end{lemma}

\medskip

\noindent
\textit{Proof of Lemma \ref{ll2}.}
Let us first notice that we can write from (\ref{10}) and (\ref{12}) that for all $X$ and $Y$ in $\rr^{2n}$,
$$q_j(X;Y)=\sigma(X,F_jY)=-\sigma(F_jX,Y)=\langle(-\sigma F_j)X,Y \rangle,$$
for $j \in \{1,2\}$ if 
\begin{equation}\label{inf13.5}\inc
\sigma=\left(
  \begin{array}{cc}
  0 & I_n \\
  -I_n & 0 \\
  \end{array}
\right) \textrm{ and } \langle X,Y \rangle=\sum_{j=1}^{2n}{X_j Y_j}.\num
\end{equation}
It follows that for all $X \in \rr^{2n}$,
\begin{equation}\label{inf14}\inc
\nabla q_j(X)=-2 \sigma F_j X. \num
\end{equation}
Since we have
$$\{q_1,q_2\}(X)=\sigma\big(\nabla q_1(X),\nabla q_2(X)\big), \ X=(x,\xi) \in \rr^{2n},$$ 
we deduce from (\ref{11}), (\ref{12}), (\ref{inf13.5}) and (\ref{inf14}) that for all $X \in \rr^{2n}$,
\inc
\begin{multline*}\label{inf15}
\{q_1,q_2\}(X)=\sigma(-2 \sigma F_1 X,-2 \sigma F_2 X)=4\langle\sigma^2 F_1 X,\sigma F_2 X \rangle
=-4 \langle F_1 X, \sigma F_2 X \rangle \\ 
=-4 \langle \sigma^t F_1 X,  F_2 X \rangle=4 \langle \sigma F_1 X,  F_2 X \rangle=4\sigma(F_1 X,F_2 X)=-4\sigma(X,F_1 F_2X),\num
\end{multline*}
which induces that 
\inc
\begin{multline*}\label{inf16}
\{q_1,q_2\}(X)=\frac{1}{2}\big[4\sigma(F_1 F_2X, X) -4\sigma(X,F_1 F_2X)\big]\\ =2\big[\sigma(X, F_2 F_1X) -\sigma(X,F_1 F_2X)\big]
=-2\sigma(X,[F_1,F_2]X), \num
\end{multline*}
by skew-symmetry of $\sigma$ and skew-symmetry of Hamilton maps with respect to $\sigma$.
Since for all $X$ and $Y$ in $\rr^{2n}$,
\begin{multline*}
\sigma(X,[F_1,F_2]Y)=\sigma(X, F_1 F_2Y) -\sigma(X,F_2 F_1Y)\\ =\sigma(F_2 F_1X, Y) -\sigma(F_1 F_2X,Y)
=-\sigma([F_1,F_2]X,Y)=\sigma(Y,[F_1,F_2]X),
\end{multline*}
we deduce that 
$$\sigma(X,-2[F_1,F_2]Y),$$ 
is the polar form associated to the quadratic form $\{q_1,q_2\}$ and that $-2[F_1,F_2]$ is its Hamilton map. $\Box$

\bigskip

Let us now consider the following quadratic form
\begin{equation}\label{2.3.98}\inc
r(X)=-\sum_{j=0}^{2n-1}{\textrm{Re } q\big((\textrm{Im } F)^j X\big)}, \ X \in \rr^{2n}. \num
\end{equation}
Since from (\ref{inf6}), $\textrm{Re } q$ is a non-positive quadratic form, we already know that $r$ is a non-negative quadratic form. We are now going to see that  
$r$ is actually a positive definite quadratic form.

\medskip

\begin{lemma}\label{l2}
The quadratic form $r$ is positive definite. 
\end{lemma}

\medskip

\noindent
\textit{Proof of Lemma \ref{l2}.} Let us assume that there exists $X_0 \in \rr^{2n}$, $X_0 \neq 0$ such that
$$r(X_0)=0.$$
The non-positivity of the quadratic form $\textrm{Re } q$ induces that for all $j=0,...,2n-1$,
\begin{equation}\label{2.3.99}\inc
\textrm{Re }q\big((\textrm{Im } F)^j X_0\big)=0. \num
\end{equation}
If $\textrm{Re }q(X;Y)$ stands for the polar form associated to the quadratic form $\textrm{Re }q$, we deduce from the Cauchy-Schwarz inequality and (\ref{2.3.99}) that
for all $j=0,...,2n-1$ and $Y \in \rr^{2n}$, 
\begin{align*}
|\textrm{Re }q\big(Y;(\textrm{Im F})^j X_0\big)|^2= & \ |\sigma\big(Y,\textrm{Re }F (\textrm{Im }F)^j X_0\big)|^2 \\
 \leq & \ [-\textrm{Re }q(Y)] [-\textrm{Re }q\big((\textrm{Im } F)^j X_0\big) ]=0. 
\end{align*} 
It follows that for all $j=0,...,2n-1$ and $Y \in \rr^{2n}$, 
$$\sigma\big(Y,\textrm{Re }F (\textrm{Im }F)^j X_0\big)=0,$$
which implies that for all $j=0,...,2n-1$, 
\begin{equation}\label{2.3.100}\inc
\textrm{Re }F(\textrm{Im }F)^j X_0=0, \num 
\end{equation}
because the symplectic form is non-degenerate. 
Let us set 
\begin{equation}\label{2.3.101}\inc
\left\lbrace
  \begin{array}{c}
  e_1=X_0\\
\displaystyle \eps_1=-\frac{1}{\textrm{Im } q(X_0)}\textrm{Im } F X_0.
  \end{array} \right. \num
\end{equation}
This is licit. Indeed, since from (\ref{2.3.99}), $\textrm{Re } q(X_0)=0$, the ellipticity of the quadratic form $q$ implies that 
$$\textrm{Im } q(X_0) \neq 0,$$
because $X_0 \neq 0$. 
By using the skew-symmetry of the Hamilton map $\textrm{Im }F$, it follows that 
$$\sigma(\eps_1,e_1)=\sigma\big{(}- (\textrm{Im } q(X_0))^{-1} \textrm{Im } F X_0,X_0\big{)}=(\textrm{Im } q(X_0))^{-1}\sigma(X_0,\textrm{Im } F X_0)=1,$$
which shows that the system $(e_1,\eps_1)$ is symplectic.

If $n=1$, we deduce that for all $X \in \rr^{2}$,
$$\textrm{Re }q(X)=\sigma(X,\textrm{Re }FX)=0,$$
because from (\ref{2.3.100}) and (\ref{2.3.101}),
$$\textrm{Re }F e_1=\textrm{Re }F \eps_1=0.$$
This contradicts (\ref{inf6}) and ends the proof of the lemma \ref{l2} in the case where $n=1$. 
We can thus assume in the following that $n \geq 2$.

In view of the proposition $21.1.3$ in \cite{hormander}, this system can be completed in a symplectic basis 
$(e_1,...,e_n,\eps_1,...,\eps_n)$ of $\rr^{2n}$.

We now check that the vector space $S$ generated by the vectors 
$e_1$ and $\eps_1$, 
\begin{equation}\label{infas1}\inc
S=\textrm{Vect}(e_1,\eps_1),\num
\end{equation} 
is not stable by the Hamilton map $\textrm{Im }F$. Indeed, if the vector space $S$ was stable by $\textrm{Im }F$, 
we could write in the symplectic coordinates associated to the symplectic basis $(e_1,...,e_n,\eps_1,...,\eps_n)$ the decomposition
\inc
\begin{multline*}\label{2.3.102}
q(x,\xi)=q(x_1,0',\xi_1,0')+ 
2q\big((x_1,0',\xi_1,0');(0,x_2,...,x_n,0,\xi_2,...,\xi_n)\big)\\ + q(0,x_2,...,x_n,0,\xi_2,...,\xi_n) 
=  q(0,x_2,...,x_n,0,\xi_2,...,\xi_n)+i\textrm{Im } q(x_1,0',\xi_1,0'), \num
\end{multline*}
because on one hand, we have
\begin{multline*}
q(x_1,0',\xi_1,0')=\sigma\big(x_1 e_1+\xi_1 \eps_1,F(x_1 e_1+\xi_1 \eps_1)\big)\\
=i\sigma\big(x_1 e_1+\xi_1 \eps_1,\textrm{Im }F(x_1 e_1+\xi_1 \eps_1)\big)
 =i\textrm{Im } q(x_1,0',\xi_1,0'),
\end{multline*}
since from (\ref{2.3.100}) and (\ref{2.3.101}), 
\begin{equation}\label{2.3.103}\inc
\textrm{Re }F e_1= \textrm{Re }F X_0=0 \ \textrm{and }  \textrm{Re }F \eps_1=-\big(\textrm{Im } q(X_0)\big)^{-1}\textrm{Re }F \textrm{Im }F X_0=0, \num
\end{equation}
and that on the other hand
\begin{eqnarray*}
&&q\big((x_1,0',\xi_1,0');(0,x_2,...,x_n,0,\xi_2,...,\xi_n)\big)\\
&=&\sigma\big(x_1 e_1+\xi_1 \eps_1,F(x_2 e_2+...+x_n e_n+\xi_2 \eps_2+...+\xi_n
\eps_n) \big)
\\
&=&-\sigma\big( F(x_1 e_1+\xi_1 \eps_1),x_2 e_2+...+x_n e_n+\xi_2 \eps_2+...+\xi_n \eps_n)\big)\\
&=&
-i\sigma\big(\textrm{Im } F(x_1 e_1+\xi_1 \eps_1),x_2 e_2+...+x_n e_n+\xi_2 \eps_2+...+\xi_n \eps_n)\big)\\
&=&0,
\end{eqnarray*}
because we assume that the space $S$ is stable by $\textrm{Im }F$. Since the quadratic form $q$ is elliptic, we deduce 
from (\ref{2.3.102}) that the quadratic form  
$$(x_1,\xi_1) \mapsto \textrm{Im } q(x_1,0',\xi_1,0'),$$
is also elliptic on $\rr^2$. 
It follows from the result proved by J. Sjöstrand (Lemma~3.1 in \cite{sjostrand}) and recalled in (\ref{12.1}) that a real-valued elliptic quadratic form is necessarily
positive definite or negative definite. Then, we deduce from (\ref{2.3.102}) and the lemma~18.6.4 in \cite{hormander} that we can find a new symplectic basis $(\tilde{e}_1,\tilde{\eps}_1)$ of the space $S$ 
such that in these new coordinates
\begin{equation}\label{2.3.104}\inc
q(x,\xi)=q(0,x_2,...,x_n,0,\xi_2,...,\xi_n)+i\lambda(x_1^2+\xi_1^2), \num
\end{equation} 
with $\lambda \in \rr^*$. We can easily check that (\ref{2.3.104}) implies that the Hamilton map $F$ associated to the quadratic form $q$ has some real eigenvalues, 
which contradicts our assumption. 
Indeed, by computing the Hamilton map $F$ from (\ref{2.3.104}) as in (\ref{2.3.88}), we obtain that  
$$F X_1=- \lambda X_1 \ \textrm{and } F X_2= \lambda X_2,$$
if $X_1=\tilde{e}_1+i  \tilde{\eps}_1$ and $X_2=\tilde{e}_1-i \tilde{\eps}_1$. We conclude that the vector space~$S$ is actually not stable by the 
Hamilton map $\textrm{Im }F$.

We can now resume the symplectic coordinates associated to the symplectic basis $(e_1,...,e_n,\eps_1,...,\eps_n)$ defined just before (\ref{infas1}). 
Since the vector space $S$ is not stable by $\textrm{Im }F$, we can find by using the following symplectically orthogonal direct sum 
$$\rr^{2n}=S \oplus^{\sigma \perp} S^{\sigma \perp},$$
and (\ref{2.3.101}), an element $X_3 \in S^{\sigma \perp}$, $X_3 \neq 0$ and some real numbers $(\lambda,\mu) \in \rr^2$ such that 
\begin{equation}\label{2.3.105}\inc
(\textrm{Im }F)^2 X_0=X_3+ \lambda e_1+\mu \eps_1. \num
\end{equation} 
It follows from (\ref{2.3.100}), (\ref{2.3.101}) and (\ref{2.3.105}) that  
\begin{equation}\label{2.3.106}\inc
\textrm{Re }F X_3=\textrm{Re }F \textrm{Im }F X_3=0, \num
\end{equation} 
because $n \geq 2$.
We can then consider the following new vectors
\begin{equation}\label{2.3.107}\inc
\left\lbrace
  \begin{array}{c}
  e_2=X_3\\
\displaystyle \eps_2=-\frac{1}{\textrm{Im } q(X_3)}\textrm{Im } F X_3,
  \end{array} \right. \num
\end{equation}
which is possible according to the fact that from (\ref{2.3.106}),  
$$\textrm{Re } q(X_3)=\sigma(X_3,\textrm{Re }F X_3)=0,$$
and to the ellipticity of $q$.
It follows that 
\begin{equation}\label{2.3.108}\inc
\sigma(\eps_1,e_1)=\sigma(\eps_2,e_2)=1, \ \sigma(e_1,e_2)=\sigma(\eps_1,e_2)=0, \num 
\end{equation}
because according to our construction $e_1, \eps_1 \in S$ and $e_2 \in S^{\sigma \perp  }$ and that we get from (\ref{2.3.101}), (\ref{2.3.107}) and (\ref{2.3.108}), 
\begin{multline*}
\sigma(\eps_2,e_1)=-\big(\textrm{Im } q(X_3) \big)^{-1} \sigma(\textrm{Im } F X_3,X_0)=\big(\textrm{Im } q(X_3) \big)^{-1} \sigma(X_3,\textrm{Im } F X_0)
\\=-\big(\textrm{Im } q(X_3) \big)^{-1} \textrm{Im } q(X_0) \sigma(e_2,\eps_1)=0
\end{multline*}
and that from (\ref{2.3.101}), (\ref{2.3.105}), (\ref{2.3.107}) and (\ref{2.3.108}),
\begin{eqnarray*}
\sigma(\eps_1,\eps_2)&=&\big(\textrm{Im } q(X_0) \textrm{Im } q(X_3) \big)^{-1} \sigma(\textrm{Im } F X_0, \textrm{Im } F X_3)\\
&=&-\big(\textrm{Im } q(X_0) \textrm{Im } q(X_3) \big)^{-1} \sigma\big((\textrm{Im } F)^2 X_0, X_3\big)\\
&=&-\big(\textrm{Im } q(X_0) \textrm{Im } q(X_3) \big)^{-1} \sigma(e_2+ \lambda e_1+\mu \eps_1, e_2)\\
&=&0.
\end{eqnarray*} 
We deduce from (\ref{2.3.100}), (\ref{2.3.101}), (\ref{2.3.106}) and (\ref{2.3.107}) that  
$$\textrm{Re }Fe_1=\textrm{Re }F e_2=\textrm{Re }F \eps_1=\textrm{Re }F \eps_2=0.$$
By checking as before that the vector space 
$$\textrm{Vect}(e_1,e_2,\eps_1,\eps_2),$$
cannot be stable by the Hamilton map $\textrm{Im }F$, because it would imply the existence of real eigenvalues for $F$, we can using (\ref{2.3.100}) iterate this construction to build a 
symplectic basis $(e_1,...,e_n,\eps_1,...,\eps_n)$ of $\rr^{2n}$ 
such that for all $j=1,...,n$,
$$\textrm{Re }Fe_j=\textrm{Re }F\eps_j=0.$$
This induces that for all $X \in \rr^{2n}$,
$$\textrm{Re }q(X)=\sigma(X,\textrm{Re }FX)=0,$$
which contradicts (\ref{inf6}) and ends the proof of the lemma \ref{l2}. $\Box$

\bigskip

By using the two previous lemmas, we can now prove that the quadratic symbol~$q$ has a \textit{finite} order $\tau$ verifying
\begin{equation}\label{stefania1}\inc 
1 \leq \tau \leq 4n-2, \num 
\end{equation}
in every non-zero point of the numerical range 
\begin{equation}\label{stefania1.5}\inc 
\Sigma(q) \setminus \{0\}=q(\rr^{2n}) \setminus \{0\}. \num 
\end{equation}

Let us first recall the classical definition of the order $k(x_0,\xi_0)$
of a symbol $p(x,\xi)$ at a point $(x_0,\xi_0) \in \rr^{2n}$ (see section 27.2, chapter 27 in \cite{hormander}). This order $k(x_0,\xi_0)$ is an element of the set $\nn \cup \{+\infty\}$ defined by
\begin{equation}\label{t10.5}\inc
k(x_0,\xi_0)=\sup\big{\{}j \in \mathbb{Z} : p_I(x_0,\xi_0)=0, \ \forall\ 1 \leq |I| \leq j\big{\}}, \num
\end{equation}
where $I=(i_1,i_2,...,i_k) \in \{1,2\}^k$, $|I|=k$ and $p_I$ stands for the iterated Poisson brackets
$$p_{I}=H_{p_{i_1}}H_{p_{i_2}}...H_{p_{i_{k-1}}}p_{i_k},$$
where $p_1$ and $p_2$ are respectively the real and the imaginary part of the symbol $p$, $p=p_1+ip_2$. The order of a symbol $q$ at a point $z$ is then defined as the 
maximal order of the symbol $p=q-z$ at every point $(x_0,\xi_0) \in \rr^{2n}$ verifying 
$$p(x_0,\xi_0)=q(x_0,\xi_0)-z=0.$$

Let $X_0$ be in $\rr^{2n}$, $X_0 \neq 0$. Let us denote by $j_0$ the lowest integer in $\{0,...,2n-1\}$ such that 
\begin{equation}\label{stefania2}\inc
\textrm{Re } q\big((\textrm{Im } F)^{j_0} X_0\big)<0, \num 
\end{equation}
given the lemma \ref{l2}. It follows that for all $0 \leq j \leq j_0-1$,
\begin{equation}\label{stefania3}\inc
\textrm{Re } q\big((\textrm{Im } F)^{j} X_0\big)=0. \num 
\end{equation}
Let us prove that 
\begin{equation}\label{stefania4}\inc
H_{\textrm{Im } q}^{2j_0} \textrm{Re } q(X_0) \neq 0. \num 
\end{equation}
Indeed, we deduce from the lemma \ref{ll2} that the Hamilton map associated to the quadratic form $H_{\textrm{Im } q}^{2j_0} \textrm{Re } q$ is
\begin{equation}\label{stefania4.1}\inc
4^{j_0}[\textrm{Im }F,[\textrm{Im }F,[...,[\textrm{Im }F,\textrm{Re }F]...], \num
\end{equation}
with exactly $2j_0$ terms $\textrm{Im }F$ appearing in the previous formula. We can write 
\inc
\begin{equation}\label{stefania4.5}
4^{j_0}[\textrm{Im }F,[\textrm{Im }F,[...,[\textrm{Im }F,\textrm{Re }F]...]=\sum_{j=0}^{2j_0}{c_j (\textrm{Im }F)^j \textrm{Re }F (\textrm{Im }F)^{2j_0-j}}, \num
\end{equation}
with $c_j \in \rr^*$ for all $j=0,...,2j_0$. Indeed, by using the following identity 
$$[P,[P,Q]]=P^2Q-2PQP+QP^2,$$
we can prove by induction that for all $n \in \nn^*$, there exist some positive constants $d_{n,j}$, $j=0,...,2n$ such that 
$$[P,[P,[...,[P,Q]...]=\sum_{j=0}^{2n}{(-1)^j d_{n,j} P^j Q P^{2n-j}},$$
if there are exactly $2n$ terms $P$ in the left-hand side of the previous identity. 
It follows from (\ref{stefania4.1}) and (\ref{stefania4.5}) that 
\inc
\begin{align*}\label{stefania5}
H_{\textrm{Im } q}^{2j_0} \textrm{Re } q(X_0)= & \ c_{j_0}\sigma\big(X_0,(\textrm{Im }F)^{j_0}\textrm{Re F}(\textrm{Im }F)^{j_0}X_0 \big) \num \\ 
+ & \ \sum_{j=0}^{j_0-1}c_j \sigma\big(X_0,(\textrm{Im }F)^j\textrm{Re F}(\textrm{Im }F)^{2j_0-j}X_0 \big)\\
+ & \ \sum_{j=0}^{j_0-1}c_{2j_0-j} \sigma\big(X_0,(\textrm{Im }F)^{2j_0-j}\textrm{Re F}(\textrm{Im }F)^{j}X_0 \big). 
\end{align*}
Since on one hand
\begin{align*}
\sigma\big(X_0,(\textrm{Im }F)^{j_0}\textrm{Re F}(\textrm{Im }F)^{j_0}X_0 \big)= & \ (-1)^{j_0}\sigma\big((\textrm{Im }F)^{j_0}X_0,\textrm{Re F}(\textrm{Im }F)^{j_0}X_0 \big)\\
= & \ (-1)^{j_0}\textrm{Re }q\big((\textrm{Im }F)^{j_0}X_0\big)
\end{align*}
by skew-symmetry of the Hamilton map $\textrm{Im }F$, and that on the other hand we get from (\ref{inf6}), (\ref{stefania3}) and the Cauchy-Schwarz inequality
\begin{align*}
& \ |\sigma\big(X_0,(\textrm{Im }F)^j\textrm{Re F}(\textrm{Im }F)^{2j_0-j}X_0 \big)|\\
=& \ |\sigma\big((\textrm{Im }F)^jX_0,\textrm{Re F}(\textrm{Im }F)^{2j_0-j}X_0 \big)|\\
= & \ -\textrm{Re }q\big((\textrm{Im }F)^jX_0;(\textrm{Im }F)^{2j_0-j}X_0\big) \\
\leq & \ [-\textrm{Re }q((\textrm{Im }F)^jX_0)]^{\frac{1}{2}} [-\textrm{Re }q((\textrm{Im }F)^{2j_0-j}X_0)]^{\frac{1}{2}}=0
\end{align*}
and 
\begin{align*}
& \ |\sigma\big(X_0,(\textrm{Im }F)^{2j_0-j}\textrm{Re F}(\textrm{Im }F)^{j}X_0 \big)|\\
=& \ |\sigma\big((\textrm{Im }F)^{2j_0-j}X_0,\textrm{Re F}(\textrm{Im }F)^{j}X_0 \big)|\\
= & \ -\textrm{Re }q\big((\textrm{Im }F)^{2j_0-j}X_0;(\textrm{Im }F)^{j}X_0\big) \\
\leq & \ [-\textrm{Re }q((\textrm{Im }F)^{2j_0-j}X_0)]^{\frac{1}{2}} [-\textrm{Re }q((\textrm{Im }F)^{j}X_0)]^{\frac{1}{2}}=0
\end{align*}
if $j=0,...,j_0-1$, we deduce (\ref{stefania4}) from (\ref{stefania2}) and (\ref{stefania5}) that 
$$H_{\textrm{Im } q}^{2j_0} \textrm{Re } q(X_0)=(-1)^{j_0}c_{j_0}\textrm{Re }q\big((\textrm{Im }F)^{j_0}X_0\big) \neq 0,$$
because $c_{j_0} \in \rr^*$. We conclude from (\ref{t10.5}), (\ref{stefania2}) and (\ref{stefania4}) that the quadratic symbol~$q$ has a \textit{finite} order lower or equal than $4n-2$ 
in every non-zero point of the numerical range $\Sigma(q) \setminus \{0\}$.

\subsubsection{Some estimates for the resolvent of elliptic quadratic differential operators}
To obtain the result of exponential decay in time for the norm of the contraction semigroup given by the theorem \ref{theorem}, we need to estimate 
the resolvent of its generator in some particular regions of the resolvent set. 
The first estimate is given by the following proposition, which shows that the resolvent of an elliptic quadratic differential operator 
cannot blow up in norm far from the numerical range of its symbol. This result is proved in \cite{karel}. For the sake of completeness of this paper, we recall its proof.

\medskip

\begin{proposition}\label{18}
Let $q : \rr^n \times \rr^n \rightarrow \cc$, $n \in \nn^*$, be a complex-valued elliptic quadratic form. We have 
$$\forall z \not\in \Sigma(q), \ \big\|\big{(}q(x,\xi)^w-z\big{)}^{-1}\big\| \leq \frac{1}{d\big{(}z,\Sigma(q)
\big{)}},$$
where $d\big{(}z,\Sigma(q)\big{)}$ stands for the distance from $z$ to the numerical range $\Sigma(q)$.
\end{proposition}

\medskip
  
\noindent
\textit{Proof of Proposition \ref{18}}.  
If the numerical range is equal to the whole complex plane, there is nothing to prove. If $\Sigma(q) \neq \cc$, we have seen after (\ref{12.1}) that 
the numerical range is necessarily a closed angular sector with a top in $0$ and an opening strictly lower than $\pi$. 
Let us consider $z \not\in \Sigma(q)$
and denote by $z_0$ its orthogonal projection on the non-empty closed convex set $\Sigma(q)$. According to the shape of the numerical range, it follows that $z_0$ 
belongs to its boundary and that we can find a complex number $z_1 \in \cc^*$, $|z_1|=1$ such that  
$$\Sigma(z_1 q) \subset \big\{z \in \cc : \textrm{Re } z \geq 0\big\}$$
and
\begin{equation}\label{37}\inc
z_1 z \in \big\{z \in \cc : \textrm{Re } z<0\big\}, \ d\big{(}z,\Sigma(q)\big{)}=d(z_1 z,i\rr). \num
\end{equation}
Using now that the operator $i [\textrm{Im}(z_1 q)]^w$ 
is formally skew-selfadjoint, we obtain that for all $u \in \mathcal{S}(\rr^n)$,
\inc
\begin{equation}\label{38}
 \textrm{Re}\big{(}z_1 q(x,\xi)^w u - z_1 z u,u \big{)}_{L^2}
= d(z_1
z,i\rr)\|u\|_{L^2}^2+\big{(}\big[\textrm{Re}\big{(}z_1q(x,\xi)\big{)}\big]^w u,u\big{)}_{L^2}. \num
\end{equation}
Then, since the quadratic form $\textrm{Re}(z_1 q)$ is non-negative, we deduce from the symplectic invariance of the Weyl quantization (Theorem 18.5.9 in \cite{hormander}) 
and the theorem 21.5.3 in \cite{hormander} that there exists a metaplectic operator $U$, which is a unitary transformation of $L^2(\rr^n)$ and an automorphism of the 
spaces $\mathcal{S}(\rr^n)$ and $B$, such that  
$$\big[\textrm{Re}\big(z_1 q(x,\xi)\big)\big]^w=U^{-1} \Big(\sum_{j=1}^{\tilde{k}}{\lambda_j(D_{x_j}^2+x_j^2)}
+\sum_{j=\tilde{k}+1}^{\tilde{k}+\tilde{l}}{x_j^2} \Big) U,$$
with $\tilde{k},\tilde{l} \in \nn$ and $\lambda_j>0$ for all $j=1,...,\tilde{k}$. By using that $U$ is a unitary operator on $L^2(\rr^{n})$, we obtain that the quantity
$$\big{(} \big[\textrm{Re}\big{(}z_1q(x,\xi)\big{)}\big]^w u,u\big{)}_{L^2}= \sum_{j=1}^{\tilde{k}}{\lambda_j \big{(}\|D_{x_j} U u\|_{L^2}^2+\|x_j U u\|_{L^2}^2 \big{)}}
+\sum_{j=\tilde{k}+1}^{\tilde{k}+\tilde{l}}{\|x_j U u\|_{L^2}^2},$$
is non-negative. Then, we can deduce from the Cauchy-Schwarz inequality, $(\ref{37})$ and $(\ref{38})$ that for all $u \in \mathcal{S}(\rr^n)$,  
$$d\big{(}z,\Sigma(q) \big{)} \|u\|_{L^2} \leq |z_1| \ \|q(x,\xi)^w u -z u\|_{L^2}.$$
Finally, using the density of the Schwartz space $\mathcal{S}(\rr^n)$ in $B$ and the fact that $|z_1|=1$, we obtain that 
$$\forall z \not\in \Sigma(q), \ \big\|\big{(}q(x,\xi)^w-z\big{)}^{-1}\big\| \leq \frac{1}{d\big{(}z,\Sigma(q)\big{)}},$$
since according to (\ref{15}),
$\sigma\big( q(x,\xi)^w\big) \subset \Sigma(q). \ \Box$ 

\medskip

The second estimate that we will use is a consequence of the semiclassical result given by the following 
proposition. This result deals with a result of absence of semiclassical pseudospectrum for elliptic quadratic differential operators at the boundary of their numerical ranges. It   
is linked to some properties of subellipticity of these operators. This result is proved in \cite{karel}. We also recall its proof for the sake of completeness of this paper.

Let us consider a complex-valued quadratic form 
$$q : \rr_x^n \times \rr_{\xi}^n \rightarrow \cc$$
and $\Delta$ an open half-line with a top in $0$ included in its numerical range
$$\Delta \subset \Sigma(q).$$ 
We need to define a notion of \textit{order} for the symbol $q(x,\xi)$ on this half-line $\Delta$. 
Since $q$ is a quadratic form, all the iterated Poisson brackets are also some quadratic forms. The property of degree two homogeneity of these Poisson brackets induces
that the symbol $q$ has the same order at every point of the open half-line~$\Delta$. This allows to define the order of the symbol $q$ on the half-line~$\Delta$  by 
defining this order with this common value. As underlined in \cite{karel}, this order can be \textit{finite} or \textit{infinite}. In the case where the Weyl symbol of 
an elliptic quadratic differential operator is of finite order on an open half-line belonging to the boundary of its numerical range, we have the following result.

\medskip

\begin{proposition}\label{27}
Let us consider a complex-valued elliptic quadratic form 
$$q : \rr_x^n \times \rr_{\xi}^n \rightarrow \cc,$$
with a numerical range $\Sigma(q)$ distinct from the whole complex plane 
$$\Sigma(q) \neq \cc,$$ 
which is a closed angular sector with a top in $0$ and a \emph{positive} opening. 
If this symbol $q(x,\xi)$ is of \emph{finite} order $k$ on an open half-line belonging to the boundary of its numerical range
$$\Delta \subset \partial \Sigma(q) \setminus \{0\},$$ 
then this order is necessarily \emph{even} and there is \emph{no} semiclassical pseudospectrum of index $k/(k+1)$ on 
$\Delta$ for the associated semiclassical operator, that is that for all $z \in \Delta$,
\begin{equation}\label{4.5}\inc
\exists C>0, \exists h_0>0, \forall \ 0<h<h_0, \ \big\|\big(q(x,h\xi)^w-z\big)^{-1}\big\|<C h^{-\frac{k}{k+1}}. \num
\end{equation}
\end{proposition}

\medskip

\noindent
\textit{Proof of Proposition \ref{27}.}
Let us consider the following symbol belonging to the $C_b^{\infty}(\rr^{2n},\cc)$ space composed of bounded complex-valued functions on $\rr^{2n}$ with all derivatives bounded
\begin{equation}\label{2.3.62}\inc 
r(x,\xi)=\frac{q(x,\xi)-z}{1+x^2+\xi^2}, \num
\end{equation}
with $z \in \Delta$. Setting 
$\tilde{\Sigma}(r)=\overline{r(\rr^{2n})},$
we can first notice that 
$$z \in \partial \Sigma(q) \setminus \{0\} \Rightarrow 0 \in \partial \tilde{\Sigma}(r).$$ 
Let us also notice that the symbol $r$ fulfills the principal-type condition in $0$. Indeed, if $(x_0,\xi_0) \in \rr^{2n}$ was such that $r(x_0,\xi_0)=0$ and $dr(x_0,\xi_0)=0$, 
we would get from (\ref{2.3.62}) that  
\begin{equation}\label{2.3.63}\inc
dq(x_0,\xi_0)=0. \num
\end{equation}
Since by using the lemma 3.1 in \cite{sjostrand} and the lemma 18.6.4 in \cite{hormander}, we can find a non-zero complex number $\alpha \in \cc^*$ and a real linear 
symplectic transformation $\chi$ of $\rr^{2n}$ such that 
\begin{equation}\label{51}\inc
\textrm{Re}(\alpha q) \circ \chi(x,\xi)=\sum_{j=1}^{n}{\lambda_j(\xi_j^2+x_j^2)},\num
\end{equation}
with $\lambda_j>0$ for all $j=1,...,n$, we deduce from (\ref{2.3.63}) that 
$$d\textrm{Re}(\alpha q)(x_0,\xi_0)=2\sum_{j=1}^{n}{\lambda_j \big{(}(\tilde{x}_0)_j dx_j+ (\tilde{\xi}_0)_j d\xi_j \big{)}} \circ \chi^{-1}=0,$$
with $(\tilde{x}_0,\tilde{\xi}_0)=\chi^{-1}(x_0,\xi_0)$. This would first imply that 
$$(\tilde{x}_0,\tilde{\xi}_0)=(0,0) \textrm{ and } (x_0,\xi_0)=(0,0),$$
because $\lambda_j>0$ for all $j=1,...,n$. 
It would then follow that $q(x_0,\xi_0)=0$ because $q$ is a quadratic form. However, since $r(x_0,\xi_0)=0$, we get from (\ref{2.3.62}) that $q(x_0,\xi_0)=z \neq 0$ because  
$$z \in \Delta \subset \partial\Sigma(q) \setminus \{0\},$$ 
which induces a contradiction. It follows that the symbol $r$ actually fulfills the principal-type condition in $0$.
Let us now notice that, since symbol $q$ is of finite order~$k$ in $z$, this induces in view of (\ref{2.3.62}) that the symbol $r$ is also of finite order $k$ in $0$.
On the other hand, we deduce from (\ref{2.3.62}) and (\ref{51}) that the set
$$\{(x,\xi) \in \rr^{2n} : r(x,\xi)=0\}=\{(x,\xi) \in \rr^{2n} : q(x,\xi)=z\},$$
is compact. Under these conditions, we can apply the theorem 1.4 in \cite{dencker}, which proves that the integer~$k$ is \textit{even} and gives the existence of positive constants $h_0$ and $C_1$ 
such that 
\begin{equation}\label{2.3.64}\inc
\forall \ 0<h<h_0, \forall u \in \mathcal{S}(\rr^n), \ \|r(x,h\xi)^w u\|_{L^2} \geq C_1 h^{\frac{k}{k+1}}\|u\|_{L^2}. \num
\end{equation}

\bigskip
\noindent
\textit{Remark.} We did not check the dynamical condition $(1.7)$ in~\cite{dencker}, because this assumption is not necessary
for the proof of Theorem~1.4. Indeed, this proof only use a part of the proof of Lemma 4.1 in \cite{dencker} (a part of the second paragraph), where this condition (1.7) is not needed.

\bigskip
\noindent
By using some results of symbolic calculus given by the theorem 18.5.4 in \cite{hormander} and (\ref{2.3.62}), we can write
\begin{equation}\label{2.3.65}\inc
r(x,h\xi)^w(1+x^2+h^2\xi^2)^w=q(x,h\xi)^w-z + h r_1(x,h\xi)^w +h^2 r_2(x,h\xi)^w, \num
\end{equation}
with 
\begin{equation}\label{2.3.66}\inc
r_1(x,\xi)=-i x \frac{\partial r}{\partial \xi}(x,\xi)+i \xi \frac{\partial r}{\partial x}(x,\xi) \num
\end{equation}
and
\begin{equation}\label{2.3.67}\inc
r_2(x,\xi)=-\frac{1}{2}\frac{\partial^2 r}{\partial x^2}(x,\xi)-\frac{1}{2}\frac{\partial^2 r}{\partial \xi^2}(x,\xi). \num
\end{equation}
We can easily check from (\ref{2.3.62}) that these functions $r_1$ and $r_2$ belong to the space $C_b^{\infty}(\rr^{2n},\cc)$, and we deduce from the Calder\'on-Vaillancourt theorem 
that there exists a positive constant $C_2$ such that for all $u \in \mathcal{S}(\rr^n)$ and $0<h \leq 1$,
\begin{equation}\label{2.3.68}\inc
\|r_1(x,h\xi)^w u\|_{L^2} \leq C_2 \|u\|_{L^2} \textrm{ and }  \|r_2(x,h\xi)^w u\|_{L^2} \leq C_2 \|u\|_{L^2}. \num
\end{equation}
It follows from (\ref{2.3.64}), (\ref{2.3.65}), (\ref{2.3.68}) and the triangular inequality that for all $u \in \mathcal{S}(\rr^n)$ and $0<h<h_0$,
\begin{align*}
& \ C_1 h^{\frac{k}{k+1}} \|(1+x^2+h^2 \xi^2)^w u\|_{L^2} \\
\leq & \ \|r(x,h\xi)^w (1+x^2+h^2 \xi^2)^w u \|_{L^2} \\
\leq & \ \|q(x,h\xi)^w u -z u \|_{L^2} + C_2 h(1+h)\|u \|_{L^2}. 
\end{align*}
Since from the Cauchy-Schwarz inequality, we have for all $u \in \mathcal{S}(\rr^n)$ and $0<h \leq 1$,
\begin{align*}
\|u \|_{L^2}^2 \leq & \ \|u \|_{L^2}^2+\|xu \|_{L^2}^2+\|hD_xu \|_{L^2}^2\\
= & \ \big((1+x^2+h^2 \xi^2)^w u,u\big)_{L^2} \\ 
\leq & \ \|(1+x^2+h^2 \xi^2)^w u \|_{L^2}\|u \|_{L^2}, 
\end{align*}
we obtain that for all $u \in \mathcal{S}(\rr^n)$ and $0<h <h_0$,
\begin{equation}\label{2.3.72}\inc
C_1 h^{\frac{k}{k+1}} \|u\|_{L^2} \leq \|q(x,h\xi)^w u -z u \|_{L^2} + C_2 h(1+h)\|u \|_{L^2}. \num
\end{equation}
Since $k \geq 1$ because $z \in \Sigma(q)$\footnote{We recall that $\Sigma(q)$ is a closed angular sector.}, we deduce from (\ref{2.3.72}) that there exist some positive constants $h'_0$ and $C_{3}$ such that 
for all $0<h<h'_0$ and $u \in \mathcal{S}(\rr^n)$,  
$$\|q(x,h\xi)^w u - z u\|_{L^2} \geq C_{3} h^{\frac{k}{k+1}}\|u\|_{L^2}.$$
Using that the Schwartz space $\mathcal{S}(\rr^n)$ is dense in $B$ and that the operator 
$$q(x,h\xi)^w+z,$$
is a Fredholm operator of index $0$, we obtain that for all $0<h<h_0'$,
$$\big\|\big(q(x,h\xi)^w-z\big)^{-1}\big\| \leq C_3^{-1}h^{-\frac{k}{k+1}},$$
which ends the proof of the proposition \ref{27}. $\Box$

\medskip

By using that $q$ is a quadratic form, we obtain from the change of variables 
$$y=h^{1/2}x,$$ 
with $h>0$, the following identity between the quantum operator 
$q(x,\xi)^w$ and its associated semiclassical operator $(q(x,h\xi)^w)_{0<h \leq 1}$,
\begin{equation}\label{6}\inc 
q(x,\xi)^w-\frac{z}{h}=\frac{1}{h}\big{(}q(y,h\eta)^w-z \big{)}, \num
\end{equation}
if $z \in \cc$. It follows that when (\ref{4.5}) is fulfilled, we get from (\ref{6}) the following estimate for the resolvent of the elliptic quadratic differential operator $q(x,\xi)^w$, 
\begin{equation}\label{stefania6}\inc
\exists C>0, \exists \eta_0 \geq 1, \forall \eta \geq \eta_0, \ \big\|\big(q(x,\xi)^w-\eta e^{i \textrm{arg}z} \big)^{-1} \big\| < C
\eta^{-\frac{1}{k+1}}. \num
\end{equation}

\subsubsection{Exponential decay in time for the norm of contraction semigroups}
Let 
$$q : \rr_x^n \times \rr_{\xi}^n \rightarrow \cc, \ n \in \nn^*,$$ 
be a complex-valued elliptic quadratic form verifying 
\begin{equation}\label{ele1}\inc
\textrm{Re }q \leq 0 \textrm{ and } \exists (x_0,\xi_0) \in \rr^{2n}, \ \textrm{Re }q(x_0,\xi_0) \neq 0.\num 
\end{equation}
We consider the nilpotent case
\begin{equation}\label{ele2}\inc
(\textrm{Re }F)^2= 0, \num 
\end{equation}
where $F$ stands for the Hamilton map associated to $q$, and we assume that this Hamilton map $F$ has no real eigenvalue. This reduction is licit according to section~\ref{s3.2.1}.
As mentioned after (\ref{12.1}), the numerical range $\Sigma(q)$ is necessarily a closed angular sector with a top in $0$ and an opening strictly lower than $\pi$ when the 
assumption (\ref{ele1}) is fulfilled. We deduce from (\ref{inf2}), (\ref{inf10}), (\ref{ele1}) and the ellipticity of $q$ that the numerical range $\Sigma(q)$ is actually a closed angular sector with a positive opening
such that one of the two half-lines $i\rr_+$ or $i\rr_-$ belongs to its boundary. We will consider in the following the case where 
\begin{equation}\label{ele3}\inc
i \rr_+ \subset \partial \Sigma(q). \num
\end{equation} 
The proof of the theorem \ref{theorem} in the second case 
$$i\rr_- \subset \partial \Sigma(q),$$
will be a straightforward adaptation of the one given when (\ref{ele3}) is fulfilled.
\begin{figure}[ht]
\caption{Numerical range $\Sigma(q)$.}
\centerline{\includegraphics[scale=0.7]{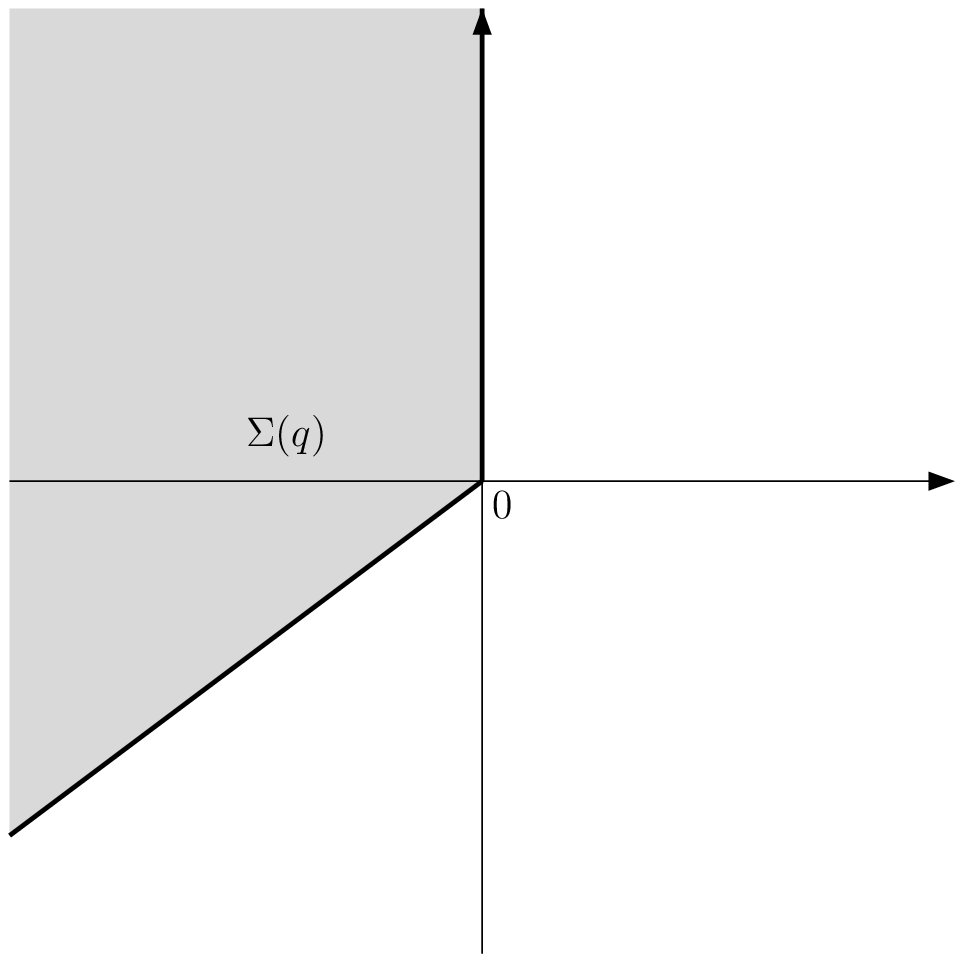}}
\end{figure}

We know from (\ref{15}) that the spectrum of the elliptic quadratic differential operator $q(x,\xi)^w$ is only composed of eigenvalues and is given by 
$$\sigma\big{(}q(x,\xi)^w\big{)}=\Big\{ \sum_{\substack{\lambda \in \sigma(F), \\  -i \lambda \in \Sigma(q) \setminus \{0\}}}
{\big{(}r_{\lambda}+2 k_{\lambda}
\big{)}(-i\lambda) : k_{\lambda} \in \nn} 
\Big\}.$$
Since from (\ref{ele1}),
$$\Sigma(q) \setminus \{0\} \subset \{z \in \cc : \textrm{Re }z \leq 0\},$$
it follows that 
$$\textrm{Re}(-i\lambda)=\textrm{Im }\lambda \leq 0,$$
if $-i\lambda \in \Sigma(q) \setminus \{0\}$. Since we assume that the Hamilton map $F$ has no real eigenvalue, we obtain that 
$$\textrm{Re}(-i\lambda)<0,$$
if $\lambda \in \sigma(F)$ and $-i\lambda \in \Sigma(q) \setminus \{0\}$. This implies that 
\begin{equation}\label{ele5}\inc
\sigma\big(q(x,\xi)^w\big) \subset \{z \in \cc : \textrm{Re }z \leq -a_0\}, \num
\end{equation}
with 
\begin{equation}\label{ele6}\inc
a_0=\inf\Big\{\sum_{\substack{\lambda \in \sigma(F), \\  -i \lambda \in \Sigma(q) \setminus \{0\}}}
{\big{(}r_{\lambda}+2 k_{\lambda}
\big{)}\big(-\textrm{Re}(-i\lambda)\big) : k_{\lambda} \in \nn} \Big\}>0. \num
\end{equation}

Let us prove that there exists a curve $\gamma \in C^{\infty}(\rr,\cc)$ in the complex plane such that there exist some positive constants $t_0$, $b$, $c_1$ and $c_2$ such that 
\begin{equation}\label{ele7}\inc
\forall t \in \rr, \ \textrm{Re } \gamma(t) \leq -b <0, \num
\end{equation}
\begin{equation}\label{ele8}\inc
\forall t \in \rr, \ |t| \geq t_0, \ \gamma(t)=-c_1 |t|+it^{4n-1}, \num
\end{equation}
\begin{equation}\label{ele9}\inc
\sigma\big(q(x,\xi)^w\big) \subset \Gamma_1 \num
\end{equation}
and 
\begin{equation}\label{ele10}\inc
\forall z \in \Gamma_2 \cup \gamma(\rr), \ \big\|\big(q(x,\xi)^w-z\big)^{-1}\big\| \leq c_2 (1+|z|)^{-\frac{1}{4n-1}}, \num
\end{equation}
where $\Gamma_1$ and $\Gamma_2$ are the subsets appearing on the following figure in the cutoff of the complex plane
\begin{equation}\label{ele11}\inc
\Gamma_1 \cup \Gamma_2 \cup \gamma(\rr)=\cc. \num
\end{equation}
\begin{figure}[ht]
\caption{Cutoff $\Gamma_1 \cup \Gamma_2 \cup \gamma(\rr)=\cc$.}
\centerline{\includegraphics[scale=0.7]{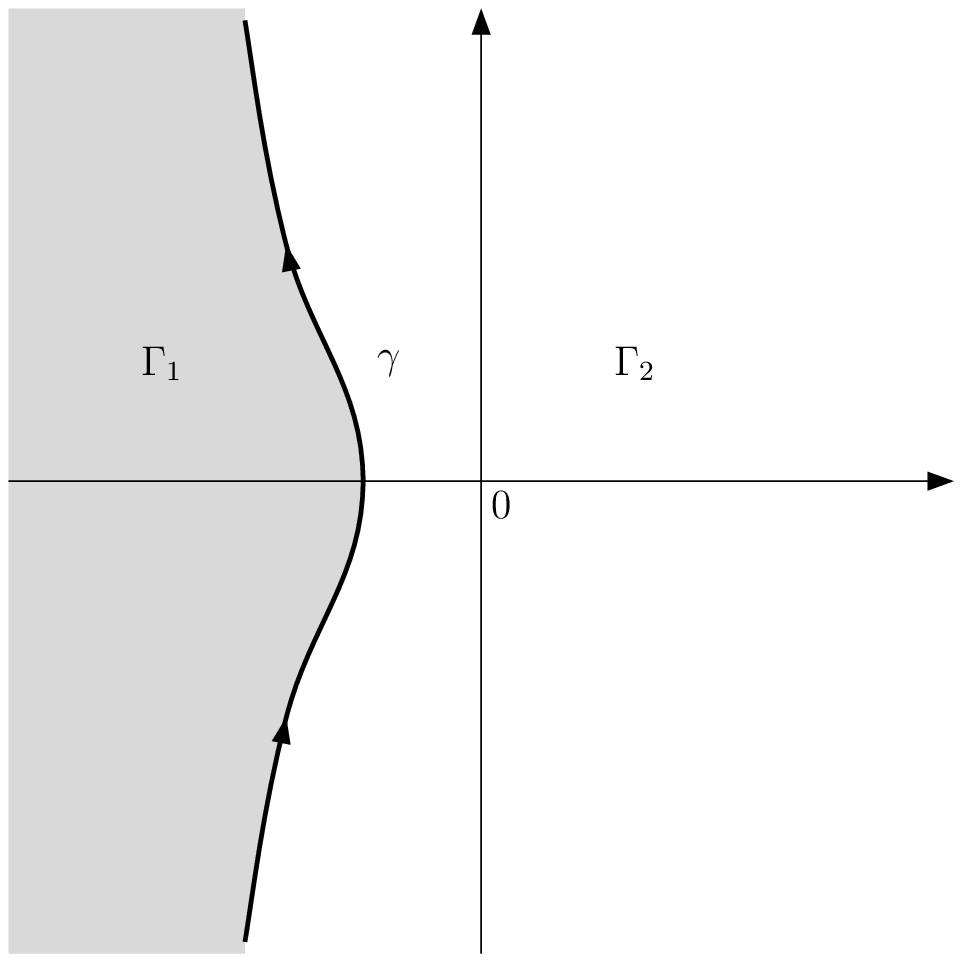}}
\end{figure}  

To justify this fact, let us begin by noticing that since we have proved in section~\ref{s3.2.2} that the quadratic symbol $q$ has always a finite order lower or equal than $4n-2$ in 
every non-zero point of the numerical range
$\Sigma(q) \setminus \{0\}$, we can deduce from the proposition \ref{27} and (\ref{stefania6}) that 
\begin{equation}\label{ele4}\inc
\exists C>0, \exists \eta_0 \geq 1, \forall \eta \geq \eta_0, \ \big\|\big(q(x,\xi)^w-i \eta \big)^{-1} \big\| < C \eta^{-\frac{1}{4n-1}}. \num
\end{equation}
This induces that for all $u \in B$ and $\eta \geq \eta_0$, 
\begin{equation}\label{ele12}\inc
\|q(x,\xi)^w u-i \eta u\|_{L^2} \geq C^{-1} \eta^{\frac{1}{4n-1}}\|u\|_{L^2}. \num
\end{equation}
It follows from (\ref{ele12}) that for all $u \in B$ and 
$$z \in \{z \in \cc : \textrm{Im }z \geq \eta_0, \ |\textrm{Re }z| \leq 2^{-1}C^{-1}|\textrm{Im }z|^{\frac{1}{4n-1}}\},$$
we have 
\begin{align*}
\|q(x,\xi)^w u -z u\|_{L^2} \geq & \ \|q(x,\xi)^w u -i \textrm{Im }z u\|_{L^2} - |\textrm{Re }z| \|u\|_{L^2} \\
\geq & \ \frac{1}{2C}|\textrm{Im }z|^{\frac{1}{4n-1}} \|u\|_{L^2},
\end{align*}
which induces that 
$$\big\|\big(q(x,\xi)^w-z \big)^{-1} \big\| \leq  2C |\textrm{Im }z|^{-\frac{1}{4n-1}},$$ 
because the spectrum of the operator $q(x,\xi)^w$ is only composed of eigenvalues. Since 
$$1+|z| \leq 1+|\textrm{Im }z|+|\textrm{Re }z| \leq 1+|\textrm{Im }z|+\frac{1}{2C}|\textrm{Im }z|^{\frac{1}{4n-1}} \leq 2|\textrm{Im }z|,$$
if the positive constant $\eta_0$ is sufficiently large and
$$z \in \{z \in \cc : \textrm{Im }z \geq \eta_0, \ |\textrm{Re }z| \leq 2^{-1}C^{-1}|\textrm{Im }z|^{\frac{1}{4n-1}}\},$$
it follows that 
\inc
\begin{multline*}\label{ele13}
\exists \eta_0>0, \forall z \in \cc, \ \textrm{Im }z \geq \eta_0, \ |\textrm{Re }z| \leq \frac{1}{2C}|\textrm{Im }z|^{\frac{1}{4n-1}}, \\
\big\|\big(q(x,\xi)^w-z \big)^{-1} \big\| \leq  2^{\frac{4n}{4n-1}}C (1+|z|)^{-\frac{1}{4n-1}}. \num
\end{multline*}
According to the shape of the numerical range (see the figure 2), we can assume that 
$$\{z \in \cc : \textrm{Im }z \leq -\eta_0, \  |\textrm{Re }z| \leq 2^{-1}C^{-1}|\textrm{Im }z|^{\frac{1}{4n-1}}  \} \subset \cc \setminus \Sigma(q),$$
if the positive constant $\eta_0$ is chosen sufficiently large. 
We can then obtain from the proposition \ref{18} that there exists a positive constant $c_3$ such that for all
$$z \in \{z \in \cc : \textrm{Im }z \leq -\eta_0, \ |\textrm{Re }z| \leq 2^{-1}C^{-1}|\textrm{Im }z|^{\frac{1}{4n-1}}\},$$ 
\begin{equation}\label{ele14}\inc
 \big\|\big(q(x,\xi)^w-z \big)^{-1} \big\| \leq  \frac{1}{d\big(z,\Sigma(q)\big)} \leq \frac{c_3}{1+|z|},\num
\end{equation}
if the positive constant $\eta_0$ is chosen sufficiently large. 
We can also deduce from the proposition \ref{18} that there exists a positive constant $c_4$ such that for all
\inc
\begin{multline*}\label{jemma1}
z \in \{z \in \cc : |\textrm{Im }z| \geq \eta_0, \ \textrm{Re }z \geq 2^{-1}C^{-1}|\textrm{Im }z|^{\frac{1}{4n-1}} \} \\
\cup \Big\{z \in \cc : |\textrm{Im }z| \leq \eta_0, \ \textrm{Re }z \geq \frac{1}{2}\Big\},\num
\end{multline*} 
\begin{equation}\label{ele14.5}\inc
 \big\|\big(q(x,\xi)^w-z \big)^{-1} \big\| \leq  \frac{1}{d\big(z,\Sigma(q)\big)} \leq c_4 (1+|z|)^{-\frac{1}{4n-1}},\num
\end{equation}
since we can find a positive constant $c_5$ such that for all $z$ belonging to the set (\ref{jemma1}),  
$$d\big(z,\Sigma(q)\big) \geq \textrm{Re }z \geq c_5(1+|z|)^{\frac{1}{4n-1}}.$$
The inclusion (\ref{ele5}) and the three estimates (\ref{ele13}), (\ref{ele14}) and (\ref{ele14.5}) allow to justify (\ref{ele7}),
(\ref{ele8}), (\ref{ele9}), (\ref{ele10}) and (\ref{ele11}).

Using this curve $\gamma$, we can now give an integral formula for the contraction semigroup generated by the operator $q(x,\xi)^w$, which allows to prove the result 
of exponential decay in time for its norm. In the following, we follow a method used by F. Hérau and F. Nier in \cite{nier}.

Let us set 
\begin{equation}\label{ele15}\inc
T(s) =\frac{1}{2i\pi} \int_{\gamma}{e^{sz}\big(z-q(x,\xi)^w\big)^{-1}dz},\num
\end{equation}
for all $s>0$. Since from (\ref{ele10}), 
\begin{equation}\label{jemma2}\inc
\forall z \in \gamma(\rr), \ \big\|e^{sz} \big(z-q(x,\xi)^w\big)^{-1}\big\| \leq c_2 e^{s \textrm{Re }z} (1+|z|)^{-\frac{1}{4n-1}}, \num
\end{equation}
we deduce from (\ref{ele8}) that 
$$T(s) \in \mathcal{L}\big(L^2(\rr^n)\big),$$
for all $s>0$. Moreover if $u_0 \in B$, we can write that for all $z \not\in \sigma\big(q(x,\xi)^w\big)$,
\begin{equation}\label{ele16}\inc
\big(z-q(x,\xi)^w\big)^{-1}u_0=\frac{1}{z-1}\big[\big(z-q(x,\xi)^w\big)^{-1}+\big(q(x,\xi)^w-1\big)^{-1} \big]\big(q(x,\xi)^w-1\big)u_0,\num
\end{equation}
because from (\ref{ele5}), $1 \not\in \sigma\big(q(x,\xi)^w\big)$. Let us set
\begin{equation}\label{ele17}\inc
u(s)=\frac{1}{2i\pi}\int_{\gamma}{e^{sz}\big(z-q(x,\xi)^w\big)^{-1}u_0 dz}, \num
\end{equation}
for all $s>0$. Since we can easily check from (\ref{ele8}) that for all $s>0$,
\begin{equation}\label{ele17.5}\inc
\frac{1}{2i\pi}\int_{\gamma}{\frac{e^{sz}}{z-1}dz}=0, \num
\end{equation}
we deduce from (\ref{ele16}) and (\ref{ele17}) that for all $s>0$,
\begin{equation}\label{ele18}\inc
u(s)=\frac{1}{2i\pi}\int_{\gamma}{\frac{e^{sz}}{z-1}\big(z-q(x,\xi)^w\big)^{-1}\big(q(x,\xi)^w-1\big)u_0 dz}. \num
\end{equation}
By using (\ref{ele8}) and (\ref{jemma2}), we can check that the right-hand side of (\ref{ele18}) is continuous with respect to $s \geq 0$ and that for all $s>0$, 
\inc
\begin{align*}\label{ele19}
\frac{du(s)}{ds}= & \ \frac{1}{2i\pi}\int_{\gamma}{e^{sz}\frac{z}{z-1}\big(z-q(x,\xi)^w\big)^{-1}\big(q(x,\xi)^w-1\big)u_0 dz}\\
= & \ \frac{1}{2i\pi}\int_{\gamma}{\frac{e^{sz}}{z-1}q(x,\xi)^w\big(z-q(x,\xi)^w\big)^{-1}\big(q(x,\xi)^w-1\big)u_0 dz}, \num
\end{align*}
because 
\begin{align*} 
& \ z\big(z-q(x,\xi)^w\big)^{-1}\big(q(x,\xi)^w-1\big)u_0\\
= & \ \big(z-q(x,\xi)^w+q(x,\xi)^w\big) \big(z-q(x,\xi)^w\big)^{-1}\big(q(x,\xi)^w-1\big)u_0\\
= & \ \big(q(x,\xi)^w-1\big)u_0+q(x,\xi)^w\big(z-q(x,\xi)^w\big)^{-1}\big(q(x,\xi)^w-1\big)u_0
\end{align*}
and that from (\ref{ele17.5}),
$$\frac{1}{2i\pi}\int_{\gamma}{\frac{e^{sz}}{z-1}dz}=0.$$
Since 
\begin{align*}
& \ \frac{1}{2i\pi}\int_{\gamma}{\frac{e^{sz}}{z-1}q(x,\xi)^w\big(z-q(x,\xi)^w\big)^{-1}\big(q(x,\xi)^w-1\big)u_0 dz} \\
= & \ q(x,\xi)^w \Big( \frac{1}{2i\pi}\int_{\gamma}{\frac{e^{sz}}{z-1}\big(z-q(x,\xi)^w\big)^{-1}\big(q(x,\xi)^w-1\big)u_0 dz} \Big),
\end{align*}
because $q(x,\xi)^w \in \mathcal{L}\big(B,L^2(\rr^n)\big)$, we deduce from (\ref{ele18}) and (\ref{ele19}) that for all $s>0$,
\begin{equation}\label{ele20}\inc
\frac{du(s)}{ds}=q(x,\xi)^w u(s), \num
\end{equation}
in $L^2(\rr^n)$. On the other hand, it follows from (\ref{ele18}) that  
\begin{equation}\label{ele21}\inc
\lim_{s \rightarrow 0^+}u(s)=\Big( \frac{1}{2i\pi}\int_{\gamma}{\big(z-q(x,\xi)^w\big)^{-1}\frac{dz}{z-1} }\Big) \big(q(x,\xi)^w-1\big)u_0.  \num
\end{equation}
We need the following lemma to conclude.

\medskip

\begin{lemma}\label{lem1}
We have
$$\frac{1}{2i\pi}\int_{\gamma}{\big(z-q(x,\xi)^w \big)^{-1}\frac{dz}{z-1}}=-\big(1-q(x,\xi)^w \big)^{-1}.$$
\end{lemma}

\medskip

By using this lemma, we get from (\ref{ele17}), (\ref{ele20}) and (\ref{ele21}) that for all $u_0 \in B$ and $s>0$,
$$e^{sq(x,\xi)^w}u_0=\frac{1}{2i\pi}\int_{\gamma}{e^{sz}\big(z-q(x,\xi)^w\big)^{-1}u_0dz}.$$ 
By using that $T(s) \in \mathcal{L}\big(L^2(\rr^n) \big)$ if $s>0$, and the density of the space $B$ in $L^2(\rr^n)$, we obtain that for all $u \in L^2(\rr^n)$ and $s>0$,
\begin{equation}\label{ele22}\inc
e^{s q(x,\xi)^w}u=\frac{1}{2i\pi}\int_{\gamma}{e^{sz}\big(z-q(x,\xi)^w\big)^{-1}u dz}.  \num
\end{equation}
Let us notice from (\ref{ele7}) and (\ref{ele10}) that for all $z \in \gamma(\rr)$, $u \in L^2(\rr^n)$ and $s>0$,
\inc
\begin{equation}\label{ele23}
e^{\frac{b}{2}s}\big\|e^{sz}\big(z-q(x,\xi)^w\big)^{-1}u \big\|_{L^2} \leq c_2 (1+|z|)^{-\frac{1}{4n-1}}e^{s(\textrm{Re }z+\frac{b}{2})}\|u\|_{L^2} \leq c_2\|u\|_{L^2}. \num
\end{equation}
We deduce from (\ref{ele8}), (\ref{ele22}) and (\ref{ele23}) that for all $u \in L^2(\rr^n)$ and $s>0$,
\inc
\begin{multline*}\label{ele24}
e^{\frac{b}{2}s}\|e^{s q(x,\xi)^w}u\|_{L^2} \leq \frac{c_2}{2\pi}\Big( \int_{-t_0}^{t_0}{|\gamma'(t)|dt} \Big)\|u\|_{L^2} \\
+\frac{c_2}{\pi}\Big(\int_{t_0}^{+\infty}{\frac{\big(c_1^2+(4n-1)^2t^{8n-4}\big)^{\frac{1}{2}}}{\big(1+(c_1^2 t^2+t^{8n-2})^{\frac{1}{2}})^{\frac{1}{4n-1}}}e^{s(-c_1 t+\frac{b}{2})}dt} \Big)\|u\|_{L^2}.\num
\end{multline*}
Since 
\begin{align*}
\ & \int_{t_0}^{+\infty}{\frac{\big(c_1^2+(4n-1)^2t^{8n-4}\big)^{\frac{1}{2}}}{\big(1+(c_1^2 t^2+t^{8n-2})^{\frac{1}{2}}\big)^{\frac{1}{4n-1}}}e^{s(-c_1 t+\frac{b}{2})}dt}\\
\leq & \ \int_{t_0}^{+\infty}{\big(c_1^2+(4n-1)^2t^{8n-4}\big)^{\frac{1}{2}}e^{s(-c_1 t+\frac{b}{2})}dt} \\
\leq & \  \int_{0}^{+\infty}{\Big(c_1^2+(4n-1)^2 c_1^{-8n+4} \Big(\frac{t}{s}+\frac{b}{2}\Big)^{8n-4}\Big)^{\frac{1}{2}}e^{-t}\frac{dt}{c_1 s}},
\end{align*}
because according to (\ref{ele7}) and (\ref{ele8}),
$$s\big(c_1 t_0 -\frac{b}{2}\big) \geq \frac{s b}{2}>0,$$
we obtain from (\ref{ele24}) that there exists a positive constant $c_6$ such that for all $s \geq 1$,
\begin{equation}\label{ele25}\inc
\|e^{s q(x,\xi)^w}\|_{\mathcal{L}(L^2)} \leq c_6 e^{-\frac{b}{2}s}. \num
\end{equation}
Finally, by using that the elliptic quadratic differential operator $q(x,\xi)^w$ generates a contraction semigroup, we deduce from (\ref{ele25}) that 
$$\exists M, a>0, \forall t \geq 0, \ \|e^{t q(x,\xi)^w}\|_{\mathcal{L}(L^2)} \leq Me^{-at}.$$
It follows that it only remains to prove the lemma \ref{lem1} to end our proof of the theorem~\ref{theorem}.

To prove the lemma \ref{lem1}, it is sufficient to prove that for all $u$ and $v$ in $L^2(\rr^n)$,
\begin{equation}\label{ele26}\inc
\frac{1}{2i\pi}\int_{\gamma}{\big((z-q(x,\xi)^w )^{-1}u,v\big)_{L^2}\frac{dz}{z-1}}=-\big((1-q(x,\xi)^w )^{-1}u,v\big)_{L^2}. \num
\end{equation}
Since the function
$$z \mapsto \big((z-q(x,\xi)^w )^{-1}u,v\big)_{L^2},$$
is holomorphic in an open neighbourhood of the set $\Gamma_2 \cup \gamma(\rr)$, the Cauchy Theorem proves that 
\begin{equation}\label{ele27}\inc
\frac{1}{2i\pi}\int_{\gamma_R}{\big((z-q(x,\xi)^w )^{-1}u,v\big)_{L^2}\frac{dz}{z-1}}=-\big((1-q(x,\xi)^w )^{-1}u,v\big)_{L^2}, \num
\end{equation}
if $R >1$ and $\gamma_R$ stands for the oriented contour appearing on the following figure
\begin{equation}\label{ele28}\inc
\gamma_R=\gamma_{R,1} \cup \gamma_{R, 2} \cup \gamma_{R, 3} \cup \gamma_{R, 4}, \num
\end{equation}
\begin{figure}[ht]
\caption{Contour $\gamma_R$.}
\centerline{\includegraphics[scale=0.7]{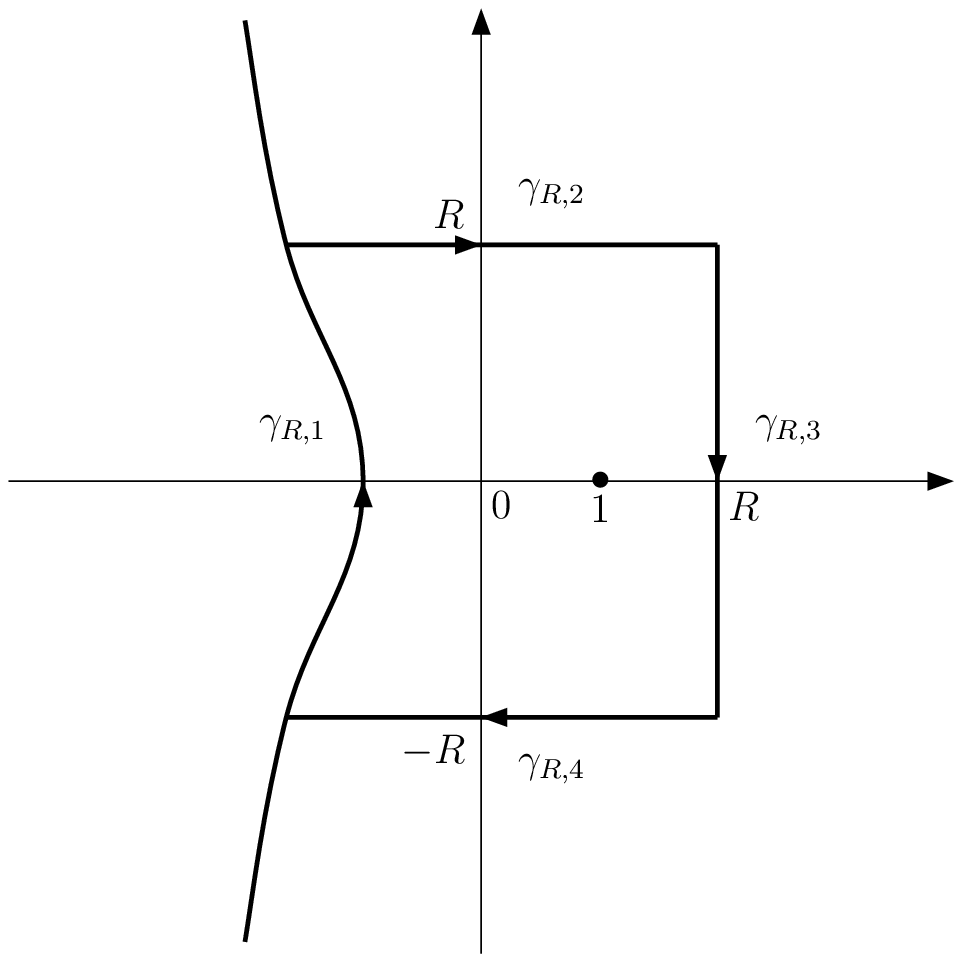}}
\end{figure}  
defined by
\begin{equation}\label{ele29}\inc
\gamma_{R,1}(t)=\gamma(t), \ -R^{\frac{1}{4n-1}} \leq t \leq R^{\frac{1}{4n-1}}, \num
\end{equation}
\begin{equation}\label{ele30}\inc
\gamma_{R,2}(t)=t+iR, \ - c_1 R^{\frac{1}{4n-1}} \leq t \leq R, \num
\end{equation}
\begin{equation}\label{ele31}\inc
\gamma_{R,3}(t)=R-it, \ -R \leq t \leq R \num
\end{equation}
and 
\begin{equation}\label{ele32}\inc
\gamma_{R,4}(t)=-t-iR, \ -R \leq t \leq c_1 R^{\frac{1}{4n-1}}. \num
\end{equation}
It follows from (\ref{ele8}) and (\ref{ele10}) that 
\begin{align*}\label{ele33}\inc
& \ \lim_{R \rightarrow +\infty} \frac{1}{2i\pi}\int_{\gamma_{R,1}}{\big((z-q(x,\xi)^w )^{-1}u,v\big)_{L^2}\frac{dz}{z-1}} \num \\
= & \ \frac{1}{2i\pi}\int_{\gamma}{\big((z-q(x,\xi)^w )^{-1}u,v\big)_{L^2}\frac{dz}{z-1}}. 
\end{align*}
For $j \geq 2$, we get from (\ref{ele10}), (\ref{ele30}), (\ref{ele31}) and (\ref{ele32}) that 
\inc
\begin{align*}\label{ele34}
& \ \Big| \frac{1}{2i\pi}\int_{\gamma_{R,j}}{\big((z-q(x,\xi)^w )^{-1}u,v\big)_{L^2}\frac{dz}{z-1}}\Big|   \num \\
\leq & \ \frac{1}{2\pi} \sup_{z \in \gamma_{R,j}}{\frac{|\big((z-q(x,\xi)^w )^{-1}u,v\big)_{L^2} |}{|z-1|}} \textrm{ Length}(\gamma_{R,j}) \\
\leq & \ \frac{c_2}{2 \pi} \max(2R, R+c_1R^{\frac{1}{4n-1}}) \Big(\sup_{z \in \gamma_{R,j}} \frac{1}{|z-1|(1+|z|)^{\frac{1}{4n-1}}}\Big)\|u\|_{L^2}\|v\|_{L^2} \rightarrow 0, 
\end{align*}
when $R \rightarrow +\infty$. Finally, we deduce (\ref{ele26}) from (\ref{ele33}) and (\ref{ele34}) by considering the limit when $R \rightarrow +\infty$ in (\ref{ele27}). This ends 
the proof of the lemma \ref{lem1} and also the proof of the theorem \ref{theorem}. $\Box$

\bigskip
\bigskip

\noindent
\textsc{Department of Mathematics, University of California, Evans Hall,
Berkeley, CA 94720, USA}\\
\textit{E-mail address:} \textbf{karel@math.berkeley.edu}

\end{document}